\documentclass{article}
\usepackage{amssymb}
\usepackage{amsfonts}

\input{tcilatex}
\begin{document}

\begin{center}
\bigskip \textbf{The local and global dynamics of a general cancer tumor
growth model with multiphase structure}

\textbf{Veli Shakhmurov}

Department of Mechanical Engineering, Okan University, Akfirat, Tuzla 34959
Istanbul, Turkey,

E-mail: veli.sahmurov@okan.edu.tr\ 

\bigskip \textbf{Rishad Shahmurov}

shahmurov@hotmail.com

University of Alabama Tuscaloosa, AL 35487

\textbf{Abstract}
\end{center}

We present a phase-space analysis of a mathematical model of tumor growth
with an immune responses. We consider mathematical analysis of the model
equations with multipoint initial condition regarding to dissipativity,
boundedness of solutions, invariance of non-negativity, local and global
stability and the basins of attractions. We derive some features of behavior
of one of three-dimensional tumor growth models with dynamics described in
terms of densities of three cells populations: tumor cells, healthy host
cells and effector immune cells. We found sufficient conditions, under which
trajectories from the positive domain of feasible multipoint initial
conditions tend to one of equilibrium points. Here, cases of the small tumor
mass equilibrium points-the healthy equilibrium point, the \textquotedblleft
death\textquotedblright\ equilibrium point have been examined. Biological
implications of our results are discussed.


\begin{center}
\textbf{1. Introduction}
\end{center}

Beginning with this article we intend to attempt to investigate the problems
of mathematical and biological approaches to modelings of cancer growth
dynamics processes and operations. It is important to take into account
\textquotedblleft the nonlinear property of cancer growth
processes\textquotedblright\ in construction of mathematical logistic
models. This nonlinearity approach appears very convenient to display
unexpected dynamics in cancer growth processes expressed in different
reactions of the dynamics to different concentrations of immune cells at
different stages of cancer growth developments $\left[ 1-21\right] $. Taking
into account all the complex processes, nonlinear mathematical models can be
estimated capable of compensation and minimization the inconsistencies
between different mathematical models related to cancer growth-anticancer
factor affections. The elaboration of mathematical non-spatial models of the
cancer tumor growth in the broad framework of tumor immune interactions
studies is one of intensively developing areas in the modern mathematical
biology, see works $[1-9]$. Of course, the development of powerful cancer
immunotherapies requires first of all an understanding of the mechanisms
governing the dynamics of tumor growth.\ One of main reasons for creation of
non-spatial dynamical models of this nature is related to the fact that they
are described by a system of ordinary differential equations, which can be
efficiently investigated by powerful methods of qualitative theory of
ordinary differential equations and dynamical systems theory. Mathematical
models for tumour growth have been extensively studied in the literature to
understand the mechanism of the disease and to predict its future behavior.
Interactions of tumour cells with other cells of the body, i.e. healthy host
cells and immune system cells are the main components of these models and
these interactions may yield different outcomes. Some important phenomena of
the tumour progression such as tumour dormancy, creeping through, and escape
from immune surveillance have been investigated by using these models.
Kuznetsov et al. $\left[ 1\right] $ proposed a model of second order,
governed by ordinary differential equations (ODEs), which includes the
effector immune cell and the tumour cell populations.They demonstrated that
even with two cell populations, these models can provide very rich dynamics
depending on the system parameters and explained some very important aspects
of the stages of cancer progression. Three equation mathematical models of
tumor growth with an immune responses were studied e.g. in $\left[ \text{4,
5, 7, 9, 10}\right] .$ For instance, Kirschner and Panetta $\left[ \text{4}%
\right] $ examined the tumour cell growth in the presence of the effector
immune cells and the cytokine IL-2 which has an essential role in the
activation and stimulation of the immune system. de Pillis and Radunskaya $%
\left[ 5\right] $ included a normal tissue cell population in this model,
performed phase space analysis and investigated the effect of chemotherapy
treatment by using optimal control theory. In $\left[ 9\right] $,
interactions between cancer cells, effector cells, and cytokines (such as
IL-2, TGF-$\beta $, IFN-$\gamma $) studid. In $\left[ 7\right] $
interactions between cancer cells, effector cells, and normal tissue cells
are ivestigated. In $\left[ 6\right] $, a four-dimensional model is
discussed which can undergo Hopf bifurcations leading to periodic orbits, a
possible route to the development of chaotic attractors (for general review
see e.g. $\left[ \text{1, 3}\right] $). \ In $\left[ 10\right] $ global
behavior of the tumour growth population dynamics was investigated. The
local stability, the chaotic behavior properties and some features of global
behavior tumour growth model of $\left( 1.1\right) $ with the classical
initial condition were studied in $\left[ 12\right] $ and $\left[ 11\right] $%
, respectively. The complex oscillations were studied in $\left[ 16\right] $%
. Moreover, the model has been also used to define optimal control problems
(see e.g. $\left[ 16-18\right] $). Note that nonlinear dynamic systems
studied e.g. in $[22-24]$. In contrast to mentioned works, here mathematical
analysis of multipoint IVP for $\left( 1.1\right) $, local and global
stability and the multiphase basins of attractions have been investigated.
We prove that all orbits are bounded and must converge to one of several
possible equilibrium points. Therefore, the long-term behavior of an orbit
is classified according to the basin of multipoint attraction in which it
starts. Here, we examine the dynamics of one cancer growth model proposed in 
$[5],$ but possessing multiphase structure, i.e. we consider the following
multipoint initial value problem (IVP) for dynamical system

\[
\dot{x}_{1}=B_{1}\left( x_{1}\right) -D_{1}\left( x_{1},x_{2}\right)
-h_{1}\left( x_{1,}x_{3}\right) , 
\]%
\begin{equation}
\dot{x}_{2}=B_{2}\left( x_{2}\right) -D_{2}\left( x_{2}\right) -h_{2}\left(
x_{1},x_{2}\right) ,  \tag{1.1}
\end{equation}%
\[
\dot{x}_{3}=B_{3}\left( x_{1},x_{3}\right) -D_{3}\left( x_{3}\right)
-h_{3}\left( x_{1},x_{3}\right) ,\text{ }t\in \left[ 0,\right. \left.
T\right) , 
\]

\begin{equation}
x_{1}\left( t_{0}\right) =x_{10}+\dsum\limits_{k=1}^{m}\alpha
_{1k}x_{1}\left( t_{k}\right) \text{, }x_{2}\left( t_{0}\right)
=x_{20}+\dsum\limits_{k=1}^{m}\alpha _{2k}x_{2}\left( t_{k}\right) \text{, }
\tag{1.2}
\end{equation}%
\[
x_{3}\left( t_{0}\right) =x_{30}+\dsum\limits_{k=1}^{m}\alpha
_{3k}x_{3}\left( t_{k}\right) \text{, }t_{0}\in \left[ 0,\right. \left.
T\right) \text{, }t_{k}\in O_{\delta }\left( t_{0}\right) , 
\]%
where $x_{1}=x_{1}\left( t\right) ,$ $x_{2}=x_{2}\left( t\right) $, $%
x_{3}=x_{3}\left( t\right) $ denote the densities of tumor cells, healthy
host cells and the effector immune cells, respectively at the moment; $t,$ $%
\alpha _{jk}$ are real numbers, $m$ is a natural number and 
\begin{equation}
O_{\delta }\left( t_{0}\right) =\left\{ t\in \mathbb{R}:\left\vert
t-t_{0}\right\vert <\delta \right\} \text{ for a }\delta >0;  \tag{1.3}
\end{equation}%
$B_{i}\left( x_{i}\right) ,$ $i=1,2$ correspond to the logistic growth of
tumor and normal health cells in the absence of any effect from immune cells
populations; $D_{1},$ $h_{1}$ are the death rates of tumor cells
respectively, with interaction of normal and immune cells; $D_{2}$ is the
natural death rate of normal health cells $x_{2}$ and $h_{2}$ is the death
rates of $x_{2}$ with interaction of tumor cells; $D_{3}$ is the natural
death rate of immune cells $x_{3}$ and $h_{3}$ is the death rates of $x_{3}$
with interaction of tumor cells; The third equation of the model describes
the change in the immune cells population with time $t.$ The first term $%
B_{3}\left( x_{1},x_{3}\right) $ of the third equation illustrates the
stimulation of the immune system by the tumor cells with tumor specific
antigens. The rate of recognition of the tumor cells by the immune system
depends on the antigenicity of the tumor cells. The model of the recognition
process is given by the rational type function which depends on the number
of tumor cells;\ $\alpha _{jk}$ are real numbers and $m$ is a natural number
such that,%
\begin{equation}
x_{j0}+\dsum\limits_{k=1}^{m}\alpha _{jk}x_{j}\left( t_{k}\right) \geq 0%
\text{, }j=1,2,3.  \tag{1.4}
\end{equation}

Note that, for $\alpha _{j1}=\alpha _{j2}=...\alpha _{jm}=0$ and 
\[
B_{1}\left( x_{1}\right) =r_{1}x_{1}\left( 1-k_{1}^{-1}x_{1}\right) ,\text{ }%
D_{1}\left( x_{1,}x_{2}\right) =a_{12}x_{1}x_{2},\text{ }h_{1}\left(
x_{1,}x_{3}\right) =a_{13}x_{1}x_{2}, 
\]%
\[
B_{2}\left( x_{2}\right) =r_{2}x_{2}\left( 1-k_{2}^{-1}x_{2}\right) ,\text{ }%
D_{2}\left( x_{2}\right) =0,\text{ }h_{2}\left( x_{1,}x_{2}\right)
=a_{21}x_{1}x_{2}, 
\]%
\[
B_{3}\left( x_{1},x_{3}\right) =\frac{r_{3}x_{1}x_{3}}{x_{1}+k_{3}},\text{ }%
D_{3}\left( x_{1,}x_{3}\right) =d_{3}x_{3},\text{ }h_{3}\left(
x_{1},x_{3}\right) =a_{31}x_{1}x_{3} 
\]%
the problem $\left( 1.3\right) -\left( 1.4\right) $ becomes the following IVP

\[
\dot{x}_{1}=r_{1}x_{1}\left( 1-k_{1}^{-1}x_{1}\right)
-a_{12}x_{1}x_{2}-a_{13}x_{1}x_{3}, 
\]%
\begin{equation}
\dot{x}_{2}=r_{2}x_{2}\left( 1-k_{2}^{-1}x_{2}\right) +a_{21}x_{1}x_{2} 
\tag{1.5}
\end{equation}%
\[
\dot{x}_{3}=\frac{r_{3}x_{1}x_{3}}{x_{1}+k_{3}}-a_{31}x_{1}x_{3}-d_{3}x_{3},%
\text{ }t\in \left[ 0,T\right] , 
\]%
\[
x_{1}\left( t_{0}\right) =x_{10}\text{, }x_{2}\left( t_{0}\right) =x_{20},%
\text{ }x_{3}\left( t_{0}\right) =x_{30},\text{ }t_{0}\in \left[ 0,\right.
\left. T\right) 
\]%
considered in $\left[ 5\right] ,$ where $a_{ij}$, $r_{i}$, $d_{3}$\ are
positive numbers, $\alpha _{jk}$ are real numbers and $m$ is a natural
number such that 
\[
x_{1}\left( t_{0}\right) >0,\text{ }x_{2}\left( t_{0}\right) >0,\text{ }%
x_{3}\left( t_{0}\right) >0, 
\]%
where the first term of the first equation corresponds to the logistic
growth of tumor cells in the absence of any effect from other cells
populations with the growth rate of $r_{1}$ and maximum carrying capacity $%
k_{1}$. The competition between host cells and tumor cells $x_{1}\left(
t\right) $ which results in the loss of the tumor cells population is given
by the term $a_{12}x_{1}x_{2}$. Next, the parameter $a_{13}$ refers to the
tumor cell killing rate by the immune cells $x_{3}\left( t\right) $. In the
second equation, the healthy tissue cells also grow logistically with the
growth rate of $r_{2}$ and maximum carrying capacity $k_{2}$. We assume that
the cancer cells proliferate faster than the healthy cells which gives $%
r_{1}>r_{2}$. The tumor cells also inactivate the healthy cells at the rate
of $a_{21}$. The third equation of the model describes the change in the
immune cells population with time $t.$ The first term of the third equation
illustrates the stimulation of the immune system by the tumor cells with
tumor specific antigens. The model of the recognition process depends on the
number of tumor cells with positive constants $r_{3}$ and $k_{3}$. The
immune cells are inactivated by the tumor cells at the rate of $a_{31}$ as
well as they die naturally at the rate $d_{3}.$

We suppose that the constant influx $s$ of the activated effector cells into
the tumor microenvironment is zero. Therein, note that, the references and
nonlinear dynamic systems studied e.g. in $[14-15]$. One of main aims is
derivation of sufficient conditions under which the possible biologically
feasible dynamics is local and globally stable, and a converges to one of
equilibrium points. Since these equilibrium points have a biological sense,
we notice that understanding limit properties of dynamics of cells
populations based on solving problems $(1.1)-\left( 1.2\right) $ may be of
an essential interest for the prediction of health conditions of a patient
without a treatment, when the data (e.g. the status of blood cells shown
above) that determines the condition of the patient are compared at various
times $t_{0},t_{1},...,t_{m}$ and correlated. Note that the local and global
stability properties of $\left( 1.1\right) $ with the classical initial
condition were studied in $\left[ 8\right] $ and $\left[ 9\right] $,
respectively. We prove that all orbits are bounded and must converge to one
of several possible equilibrium points.

\begin{center}
\textbf{2}.\textbf{\ Notations and background.}
\end{center}

Consider the multipoint IVP for nonlinear equation%
\begin{equation}
\frac{du}{dt}=f\left( u\right) ,\text{ }t\in \left[ 0,T\right] ,  \tag{2.1}
\end{equation}%
\[
u\left( t_{0}\right) =u_{0}+\dsum\limits_{k=1}^{m}\alpha _{k}u\left(
t_{k}\right) \text{, }t_{0}\in \left[ 0,\right. \left. T\right) \text{, }%
t_{k}\in \left( 0,T\right) ,\text{ }t_{k}>t_{0} 
\]%
in a Banach space $X$, where $\alpha _{k}$ are complex numbers, $m$ is a
natural number and $u=u\left( t\right) $ is a $X-$valued function. Note
that, for $\alpha _{1}=\alpha _{2}=...\alpha _{m}=0$ the problem $\left(
2.1\right) $ becomes the following local Cauchy problem%
\begin{equation}
\frac{du}{dt}=f\left( u\right) ,\text{ }u\left( t_{0}\right) =u_{0},\text{ }%
t\in \left[ 0,T\right] ,\text{ }t_{0}\in \left[ 0,\right. \left. T\right) 
\text{.}  \tag{2.2}
\end{equation}

For $u_{0}\in X$ let $\bar{B}_{r}\left( u_{0}\right) $ denotes a closed ball
in $X$ with radius $r$ centered at $u_{0}$, i.e., 
\[
\bar{B}_{r}\left( u_{0}\right) =\left\{ u\in X:\left\Vert u-u_{0}\right\Vert
_{X}\leq r\right\} . 
\]

We can generalized classical Picard existence theorem for nonlinear
multipoint IVP $\left( 2.1\right) $.

By reasoning as a classical case we obtain

\textbf{Theorem 2.1. }Let $X$ be a Banach space. Suppose $f:X\rightarrow X$
satisfies local Lipschitz condition on $\bar{B}_{r}(\upsilon _{0})\subset $ $%
X$, i.e.%
\[
\left\Vert f\left( u\right) -f\left( \upsilon \right) \right\Vert _{X}\leq
L\left\Vert u-\upsilon \right\Vert _{X} 
\]%
for each $u$, $\upsilon \in \bar{B}_{r}(\upsilon _{0})$ and there exists $%
\delta >0$ such that 
\[
t_{k}\in O_{\delta }\left( t_{0}\right) =\left\{ t\in \mathbb{R}:\left\vert
t-t_{0}\right\vert <\delta \right\} , 
\]%
where 
\[
\upsilon _{0}=u_{0}+\dsum\limits_{k=1}^{m}\alpha _{k}u\left( t_{k}\right) . 
\]

Moreover, let 
\[
M=\sup\limits_{u\in \bar{B}_{r}(\upsilon _{0})}\left\Vert f\left( u\right)
\right\Vert _{X}<\infty . 
\]

Then,\ problem $\left( 2.1\right) $ has a unique continuously differentiable
local solution $u(t)$ for $t\in O_{\delta }\left( t_{0}\right) $, where $%
\delta \leq \frac{r}{M}.$

\textbf{Proof.} We rewrite the initial value problem $\left( 2.1\right) $ as
an integral equation 
\[
u=\upsilon _{0}+\dint\limits_{t_{0}}^{t}f\left( u\left( s\right) \right)
ds.\ 
\]

For $0<\eta <\frac{r}{M}$ we define the space 
\[
Y=C\left( \left[ -\eta ,\eta \right] ;\bar{B}_{r}(\upsilon _{0})\right) . 
\]

Let%
\[
Qu=\upsilon _{0}+\dint\limits_{t_{0}}^{t}f\left( u\left( s\right) \right)
ds. 
\]
First, note that if $u\in Y$ then%
\[
\left\Vert Qu-\upsilon _{0}\right\Vert _{X}\leq \left\Vert
\dint\limits_{t_{0}}^{t}f\left( u\left( s\right) \right) ds\right\Vert
_{X}\leq M\eta <r. 
\]

Hence, $Qu\in Y$ so that $Q:Y\rightarrow Y.$ Moreover, for all $u$, $%
\upsilon \in Y$\ we have 
\[
\left\Vert Qu-Q\upsilon \right\Vert _{X}\leq \left\Vert
\dint\limits_{t_{0}}^{t}\left[ f\left( u\left( s\right) \right) -f\left(
\upsilon \left( s\right) \right) \right] ds\right\Vert _{X}\leq 
\]%
\begin{equation}
L_{f}\eta \left\Vert u-\upsilon \right\Vert _{X},  \tag{2.3}
\end{equation}%
where $L_{f}$ is a Lipschitz constant for $f$ on $\bar{B}_{r}(\upsilon _{0})$%
. Hence, if we choose 
\begin{equation}
\eta <\min \left\{ \frac{r}{M},\frac{1}{L_{f}}\right\}  \tag{2.4}
\end{equation}%
then $Q$ is a contraction on $Y$ and it has a unique fixed point. Since $%
\eta $ depends only on the Lipschitz constant of $f$ and on the distance $r$
of the initial data from the boundary of $\bar{B}_{r}(\upsilon _{0})$. Then
repeated application of this result gives a unique local solution defined
for $\left\vert t-t_{0}\right\vert <\frac{r}{M}.$

\textbf{Theorem 2.2. }Let $X$ be a Banach space. Suppose that $%
f:X\rightarrow X$ satisfies global Lipschitz condition, i.e.%
\[
\left\Vert f\left( u\right) -f\left( \upsilon \right) \right\Vert _{X}\leq
L\left\Vert u-\upsilon \right\Vert _{X} 
\]%
for each $u$, $\upsilon \in X$. Moreover, let 
\[
M=\sup\limits_{u\in X}\left\Vert f\left( u\right) \right\Vert _{X}<\infty . 
\]

Then\ problem $\left( 2.1\right) $ has a unique continuously differentiable
global solution $u(t)$ for all $t\in \left[ t_{0},T\right] .$

\textbf{Proof. }The key point of proof is to show that the constant $\delta $
of Theorem 2.1 can be made independent of the $\upsilon _{0}.$ It is not
hard to see that the independence of $\upsilon _{0}$ comes through the
constant $M$ in therm $\frac{r}{M}$ in $\left( 2.4\right) $. Since in the
current case the Lipschitz condition holds globally, one can choose $r$
arbitrary large. Therefore, for any finite $M$, we can choose $r$ large
enough and by using $\left( 2.3\right) ,$ $\left( 2.4\right) $ we obtain the
assertion.

Let $X$ be a Banach space. $w\in X$ is called a critical point (or
equilibria point) for the equation $\left( 2.1\right) $ if $f\left( w\right)
=0.$

We denote the solution of the problem $\left( 2.1\right) $ by%
\[
\phi \left( t,u_{0}\right) =\phi \left( t,u\left( t_{0}\right) ,u\left(
t_{1}\right) ,...,u\left( t_{m}\right) \right) . 
\]

\bigskip \textbf{Definition 2.1. }Let $u_{0}\in X$, $u\left( t\right) =\phi
\left( t,u_{0}\right) $ be a solution of $\left( 2.1\right) $ and $w\in X$
be a critical point of $\left( 2.1\right) .$ If there exists a neighbourhood 
$O\left( w\right) \subset X$ of $w$ such that $\lim\limits_{t\rightarrow
\infty }u\left( t\right) =w$ for $u_{0}+\dsum\limits_{k=1}^{m}\alpha
_{k}u\left( t_{k}\right) \subset O\left( w\right) $, $t_{0}\in \left[
0,\right. \left. T\right) $, $t_{k}\in O_{\delta }\left( t_{0}\right) $ and
a $\delta >0$, then $w$ is called a positive multiphase attractor$.$

\textbf{Definition 2.2. }Assume $w\in X$ is a multiphase attractor point of $%
\left( 2.1\right) $ and $u\left( t\right) =\phi \left( t,u_{0}\right) $ is a
solution of $\left( 2.1\right) .$ A set $\left\{ u\text{: }%
u=u_{0}+\dsum\limits_{k=1}^{m}\alpha _{k}u\left( t_{k}\right) \right\}
\subset X$ \ is called a domain of multiphase basin (multiphase attractor or
domain of multiphase asymptotic stability) of $w$ if $\lim\limits_{t%
\rightarrow \infty }u\left( t\right) =w.$

\begin{center}
\bigskip \textbf{3. Boundedness, invariance of non-negativity, and
dissipativity}
\end{center}

In this section, we shall show that the model equation are bounded with
negative divergence, positively invariant with respect to a region in $%
R_{+}^{3}$ and dissipative. As we are interested in biologically relevant
solutions of the system, the next two results show that the positive octant
is invariant and that all trajectories in this octant are recurrent. Let 
\begin{equation}
O_{K}=\left\{ x=\left( x_{1},x_{2},x_{3}\right) \in R_{+}^{3}\text{: }0\leq
x_{i}\leq K_{i}\text{, }i=1,\text{ }2,\text{ }3\right\} ,  \tag{3.1}
\end{equation}%
where 
\[
K_{i}=\max \left\{ 1,\text{ }x_{i0}+\dsum\limits_{k=1}^{m}\alpha
_{ik}x_{1}\left( t_{k}\right) \right\} \text{, }t_{k}\in O_{\delta }\left(
t_{0}\right) ,\text{ }i=1,\text{ }2,\text{ }3. 
\]%
Consider the problem $(1.3)-\left( 1.4\right) $ with $t_{0}=0.$

\textbf{Condition 3.1. }Assume:

(1) $B_{i}\left( x_{i}\right) >0$, $D_{1}\left( x_{1},x_{2}\right) >0,$ $%
D_{2}\left( x_{2}\right) >0$, $B_{1}\left( 0\right) =D_{1}\left(
0,x_{2}\right) =0,$ $\frac{d}{dx_{1}}B_{i}\left( x_{i}\right) >0,$ $\frac{d}{%
dx_{1}}D_{1}\left( x_{1},x_{2}\right) >0,$ $\frac{d}{dx_{2}}D_{2}\left(
x_{2}\right) >0$ for $x_{i}>0,$ $i=1.2;$ moreover, $\frac{d}{dx_{1}}%
B_{1}\left( 0\right) >\frac{\partial }{\partial x_{1}}D_{1}\left(
0,x_{2}\right) $ and $\frac{d}{dx_{2}}B_{2}\left( 0\right) >\frac{d}{dx_{2}}%
D_{2}\left( 0\right) $;

(2) $h_{k}\left( x_{1},x_{3}\right) $ $>0,$ $h_{k}\left( 0,x_{3}\right) =0,$ 
$h_{k}\left( x_{1},0\right) =0,$ $h_{j}\in C^{1}\left( R_{+}^{2}\right) ,$ $%
\frac{\partial h_{k}}{\partial x_{k}}\geq 0$, $\frac{\partial h_{2}}{%
\partial x_{2}}\geq 0$ for $k=1,3$ and $x\in R_{+}^{3};$

(3) $h_{2}\left( x_{1},x_{2}\right) $ $>0,$ $h_{2}\left( x_{1},0\right) =0,$ 
$h_{2}\left( 0,x_{2}\right) =0,$ 
\[
\frac{\partial }{\partial x_{1}}h_{2}\left( 0,x_{2}\right) \neq 0,\frac{%
\partial }{\partial x_{2}}h_{2}\left( 0,x_{2}\right) =0,\text{ }\frac{%
\partial }{\partial x_{1}}h_{k}\left( 0,x_{3}\right) \neq 0, 
\]%
\[
\text{ }\frac{\partial }{\partial x_{3}}h_{k}\left( 0,x_{3}\right) =0,k=1,3%
\text{ for }x\in R_{+}^{3}; 
\]

(4) $0<B_{3}\left( x_{1},x_{3}\right) \in C^{1}\left( R_{+}^{2}\right) ,$ $%
\frac{\partial }{\partial x_{1}}B_{3}\left( x_{1},x_{3}\right) >0$, $\frac{%
\partial }{\partial x_{3}}B_{3}\left( x_{1},x_{3}\right) >0$, $B_{3}\left(
x_{1},0\right) =0,$ $B_{3}\left( 0,x_{3}\right) =0$ and $\frac{\partial }{%
\partial x_{3}}B_{3}\left( x_{1},x_{3}\right) <\frac{d}{dx_{3}}\left[
D_{3}\left( x_{3}\right) -h_{3}\left( x_{1},x_{3}\right) \right] $ for $%
x_{1},$ $x_{3}>0;$

(5) $D_{3}\left( x_{3}\right) >0,$ $D_{3}\left( 0\right) =0,$ $D_{3}\left(
.\right) \in C^{1}\left( R_{+}\right) $ and $\frac{\partial }{\partial x_{3}}%
D_{3}\left( x_{3}\right) >0$ for $x_{3}>0;$

(6) there exist constants $K_{i}>0$ such that $B_{1}\left( K_{1}\right)
=D_{1}\left( K_{1},x_{2}\right) ,$ $\frac{d}{dx_{1}}B_{1}\left( K_{1}\right)
<\frac{\partial }{\partial x_{1}}D_{1}\left( K_{1},x_{2}\right) $, $%
B_{2}\left( K_{2}\right) =D_{2}\left( K_{2}\right) $ and%
\[
\frac{d}{dx_{2}}B_{2}\left( K_{2}\right) <\frac{d}{dx_{2}}D_{2}\left(
K_{2}\right) ,\frac{\partial }{\partial x_{1}}B_{1}\left( x_{1}\right) <%
\frac{\partial }{\partial x_{1}}\left[ D_{1}\left( x_{1},x_{2}\right)
-h_{1}\left( x_{1,}x_{3}\right) \right] , 
\]
\[
\frac{d}{dx_{2}}B_{2}\left( x_{2}\right) <\frac{d}{dx_{2}}\left[ D_{2}\left(
x_{2}\right) -h_{2}\left( x_{1},x_{2}\right) \right] \text{ for }x\in
R_{+}^{3}. 
\]

\textbf{Theorem 3.1. }Let the Condition 3.1 hold.Then: (1) $O_{K}$ is
positively invariant with respect to $(1.1)-\left( 1.2\right) ;$ (2) all
solutions of the problem $(1.1)-\left( 1.2\right) $ are uniformly bounded
and are attracted into the region $O_{K}$; (3) the system $(1.1)$ is
dissipative.

\textbf{Proof. }By Theorem 2.1 there exists a unique solution of multipoint
problem $\left( 1.1\right) -\left( 1.2\right) .$

\textbf{\ (}1\textbf{) }Consider the first equation of the system $\left(
1.3\right) $:\textbf{\ \ }%
\[
\dot{x}_{1}=B_{1}\left( x_{1}\right) -D_{1}\left( x_{1},x_{2}\right)
-h_{1}\left( x_{1,}x_{3}\right) 
\]

By assumption $h_{1}\left( x_{1,}x_{3}\right) >0$\ we get 
\[
\dot{x}_{1}<B_{1}\left( x_{1}\right) -D_{1}\left( x_{1},x_{2}\right) . 
\]

But there exists $K_{1}$ such that $B_{1}\left( K_{1}\right) =D_{1}\left(
K_{1},x_{2}\right) $ for $x_{2}>0$ by hypothesis (2). Then $\dot{x}_{1}<0$
in around of $K_{1}$. Thus 
\[
x_{1}\left( t\right) \leq \max \left\{ K_{1},\text{ }x_{10}+\dsum%
\limits_{k=1}^{m}\alpha _{1k}x_{1}\left( t_{k}\right) \right\} ,\text{ }\dot{%
x}_{1}<0\text{ for }x_{1}>1. 
\]

Hence, 
\begin{equation}
\limsup\limits_{t\rightarrow \infty }x_{1}\left( t\right) \leq K_{1}. 
\tag{3.2 }
\end{equation}%
For 
\[
\dot{x}_{2}=B_{2}\left( x_{2}\right) -D_{2}\left( x_{2}\right) -h_{2}\left(
x_{1},x_{2}\right) 
\]%
a similar analysis by assumpt\i ons (1)-(4) gives 
\[
x_{2}\left( t\right) \leq \max \left\{ K_{2},\text{ }x_{20}+\dsum%
\limits_{k=1}^{m}\alpha _{2k}x_{2}\left( t_{k}\right) \right\} , 
\]%
\begin{equation}
\limsup\limits_{t\rightarrow \infty }x_{2}\left( t\right) \leq K_{2}. 
\tag{3.3 }
\end{equation}

\bigskip Now consider 
\[
\dot{x}_{3}=B_{3}\left( x_{1},x_{3}\right) -D_{3}\left( x_{3}\right)
-h_{3}\left( x_{1},x_{3}\right) . 
\]%
From $\left( 3.1\right) $ by assumpt\i ons (5) and (6) we have

\[
\dot{x}_{3}<B_{3}\left( x_{1},x_{3}\right) -D_{3}\left( x_{3}\right) <0. 
\]%
Then by reasoning as the case of $x_{1}$ we deduced 
\[
x_{3}\left( t\right) \leq \max \left\{ K_{3},\text{ }x_{30}+\dsum%
\limits_{k=1}^{m}\alpha _{1k}x_{3}\left( t_{k}\right) \right\} , 
\]%
\begin{equation}
\limsup\limits_{t\rightarrow \infty }x_{3}\left( t\right) \leq K_{3}. 
\tag{3.4}
\end{equation}%
Hence, from $\left( 3.2\right) -\left( 3.4\right) $ we obtain (1) and (2)
assertions. Now, let us show (3). Let $f_{1},$ $f_{2},$ $f_{3}$ denote the
right sides of the system $\left( 1.1\right) .$ Since 
\[
\frac{\partial f_{1}}{\partial x_{1}}+\frac{\partial f_{2}}{\partial x_{2}}+%
\frac{\partial f_{3}}{\partial x_{3}}=\frac{\partial }{\partial x_{1}}%
B_{1}\left( x_{1}\right) -\frac{\partial }{\partial x_{1}}D_{1}\left(
x_{1},x_{2}\right) -\frac{\partial }{\partial x_{1}}h_{1}\left(
x_{1,}x_{3}\right) + 
\]%
\[
\frac{d}{dx_{2}}B_{2}\left( x_{2}\right) -\frac{d}{dx_{2}}D_{2}\left(
x_{2}\right) -\frac{\partial }{\partial x_{2}}h_{2}\left( x_{1},x_{2}\right)
+ 
\]%
\[
\frac{\partial }{\partial x_{3}}B_{3}\left( x_{1},x_{3}\right) -\frac{d}{%
dx_{3}}D_{3}\left( x_{3}\right) -\frac{\partial }{\partial x_{3}}h_{3}\left(
x_{1},x_{3}\right) 
\]%
by assumpt\i ons (1)-(6) we obtain 
\[
\frac{\partial f_{1}}{\partial x_{1}}+\frac{\partial f_{2}}{\partial x_{2}}+%
\frac{\partial f_{3}}{\partial x_{3}}<0\text{ for }x\in O_{K}, 
\]%
i.e. the system $(1.1)$ is dissipative.

\begin{center}
\textbf{4. The equilibria points}
\end{center}

In this section we find the equilibria points of the system $(1.1).$ The
equilibria points of $(1.1)$ are obtained by solving the system of
corresponding isocline equations

\[
B_{1}\left( x_{1}\right) -D_{1}\left( x_{1},x_{2}\right) -h_{1}\left(
x_{1,}x_{3}\right) =0, 
\]%
\begin{equation}
B_{2}\left( x_{2}\right) -D_{2}\left( x_{2}\right) -h_{2}\left(
x_{1},x_{2}\right) =0,  \tag{4.1}
\end{equation}%
\[
B_{3}\left( x_{1},x_{3}\right) -D_{3}\left( x_{3}\right) -h_{3}\left(
x_{1},x_{3}\right) =0. 
\]%
Since we are interested in biologically relevant solutions of $\left(
4.1\right) ,$ we find sufficient conditions under which this system have
positive solutions.

\textbf{Lemma 4.1. }Assume the assumptions (1)-(5) of the condition 3.1 are
satisfied. Then 
\[
E_{1}\left( 0,0,0\right) \text{, }E_{2}\left( \bar{x}_{1},0,0\right) \text{, 
}E_{3}\left( 0,\bar{x}_{2},0\right) \text{, }E_{4}\left( \bar{x}_{1},0,\bar{x%
}_{3}\right) ,\text{ }E_{5}\left( \bar{x}_{1},\bar{x}_{2},0\right) , 
\]%
\begin{equation}
\text{ }E_{6}\left( 0,\bar{x}_{2},\bar{x}_{3}\right) \text{ }  \tag{4.2}
\end{equation}%
\ are the equilibria points, where $\bar{x}_{1},\bar{x}_{2},\bar{x}_{3}$
will be defined in bellow.

\textbf{Proof. } By assumption (4), $E_{1}$, $E_{2}$ and $E_{3}$ are
equilibria points, where $\bar{x}_{1},$ $\bar{x}_{2}$ are solutions of the
equations, respectively 
\begin{equation}
B_{1}\left( x_{1}\right) =D_{1}\left( x_{1},0\right) ,\text{ }B_{2}\left(
x_{2}\right) =D_{2}\left( x_{2}\right) .  \tag{4.3}
\end{equation}%
It remains to find the points%
\[
E_{4}\left( \bar{x}_{1},0,\bar{x}_{3}\right) ,\text{ }E_{5}\left( \bar{x}%
_{1},\bar{x}_{2},0\right) \text{, }E_{6}\left( 0,\bar{x}_{2},\bar{x}%
_{3}\right) . 
\]%
Consider the point $E_{4}\left( \bar{x}_{1},0,\bar{x}_{3}\right) ,$ i.e. $%
x_{2}=0.$ Then, by assumption (4), we get that $E_{4}\left( \bar{x}_{1},0,%
\bar{x}_{3}\right) $ is equilibria point, when $\bar{x}_{1},$ $\bar{x}_{3}$
are solution of the following system of equations 
\begin{equation}
B_{1}\left( x_{1}\right) -D_{1}\left( x_{1},0\right) -h_{1}\left(
x_{1,}x_{3}\right) =0,  \tag{4.4}
\end{equation}%
\[
B_{3}\left( x_{1},x_{3}\right) -D_{3}\left( x_{3}\right) -h_{3}\left(
x_{1},x_{3}\right) =0. 
\]

\bigskip Consider the point $E_{5}\left( \bar{x}_{1},\bar{x}_{2},0\right) ,$
i.e. $x_{3}=0.$ Then, by assumption (4), we get that $E_{5}\left( \bar{x}%
_{1},\bar{x}_{2},0\right) $ is equilibria point, when $\bar{x}_{1},$ $\bar{x}%
_{2}$ are solution of the following system of equations 
\begin{equation}
B_{1}\left( x_{1}\right) -D_{1}\left( x_{1},x_{2}\right) -h_{1}\left(
x_{1,}0\right) =0,  \tag{4.5}
\end{equation}%
\[
B_{2}\left( x_{2}\right) -D_{2}\left( x_{2}\right) -h_{2}\left(
x_{1},x_{2}\right) =0. 
\]%
The point $E_{6}\left( 0,\bar{x}_{2},\bar{x}_{3}\right) $ is equilibria
point if $\bar{x}_{2},\bar{x}_{3}$\ are solution of the system%
\begin{equation}
B_{2}\left( x_{2}\right) -D_{2}\left( x_{2}\right) -h_{2}\left(
x_{1},x_{2}\right) =0,  \tag{4.6}
\end{equation}%
\[
B_{3}\left( x_{1},x_{3}\right) -D_{3}\left( x_{3}\right) -h_{3}\left(
x_{1},x_{3}\right) =0. 
\]
Let

\[
R_{+}^{3}=\left\{ x=\left( x_{1},x_{2},x_{3}\right) \in R^{3}\text{: }x_{i}>0%
\text{, }i=1,2,3\right\} . 
\]

We now discuss the local linearized stability of the system $\left(
1.1\right) -\left( 1.2\right) $ restricted to neighborhood of the
equilibrium points $\left( 4.2\right) $. The linearized matrix of $\left(
1.1\right) $ about an arbitrary equilibrium point $E\left(
x_{1},x_{2},x_{3}\right) $ is given by 
\begin{equation}
A_{E\left( x_{1},x_{2},x_{3}\right) }=  \tag{4.7}
\end{equation}%
\[
\left[ 
\begin{array}{ccc}
\frac{dB_{1}}{dx_{1}}-\frac{\partial D_{1}}{\partial x_{1}}-\frac{\partial
h_{1}}{\partial x_{1}} & -\frac{\partial D_{1}}{\partial x_{2}} & -\frac{%
\partial h_{1}}{\partial x_{3}} \\ 
-\frac{\partial h_{2}}{\partial x_{1}} & \frac{dB_{2}}{dx_{2}}-\frac{dD_{2}}{%
dx_{2}}-\frac{\partial h_{2}}{\partial x_{2}} & 0 \\ 
\frac{\partial B_{3}}{\partial x_{1}}-\frac{\partial h_{3}}{\partial x_{1}}
& 0 & \frac{dB_{3}}{dx_{3}}-\frac{dD_{3}}{dx_{3}}-\frac{\partial h_{3}}{%
\partial x_{3}}%
\end{array}%
\right] . 
\]

By assumption (4), the linearized matrices for equilibria points $\left(
4.2\right) $\ will be correspondingly as:%
\[
A_{1}=\left[ 
\begin{array}{ccc}
a_{11} & a_{12} & 0 \\ 
0 & a_{22} & 0 \\ 
a_{31} & 0 & a_{33}%
\end{array}%
\right] ,\text{ }A_{2}=\left[ 
\begin{array}{ccc}
b_{11} & b_{12} & 0 \\ 
b_{21} & b_{22} & 0 \\ 
b_{31} & 0 & b_{33}%
\end{array}%
\right] ,\text{ }A_{3}=\left[ 
\begin{array}{ccc}
c_{11} & c_{12} & 0 \\ 
c_{21} & c_{22} & 0 \\ 
c_{31} & 0 & c_{33}%
\end{array}%
\right] 
\]

\[
A_{4}=\left[ 
\begin{array}{ccc}
d_{11} & d_{12} & d_{13} \\ 
d_{21} & d_{22} & 0 \\ 
d_{31} & 0 & d_{33}%
\end{array}%
\right] ,\text{ }A_{5}=\left[ 
\begin{array}{ccc}
k_{11} & k_{12} & 0 \\ 
k_{21} & k_{22} & 0 \\ 
k_{31} & 0 & k_{33}%
\end{array}%
\right] ,\text{ }A_{6}=\left[ 
\begin{array}{ccc}
l_{11} & l_{12} & 0 \\ 
l_{21} & l_{22} & 0 \\ 
l_{31} & 0 & l_{33}%
\end{array}%
\right] , 
\]%
where 
\[
a_{11}=\frac{\partial }{\partial x_{1}}\left[ B_{1}-D_{1}\right] \left(
0\right) -\frac{\partial h_{1}}{\partial x_{1}}\left( 0\right) \text{, }%
a_{12}=-\frac{\partial D_{1}}{\partial x_{1}}\left( 0\right) , 
\]%
\begin{equation}
\text{ }a_{22}=\frac{d}{dx_{2}}\left[ B_{2}-D_{2}\right] \left( 0\right) ,%
\text{ }a_{31}=\frac{\partial B_{3}}{\partial x_{1}}\left( 0\right) -\frac{%
\partial h_{1}}{\partial x_{1}}\left( 0\right) \text{,}  \tag{4.8}
\end{equation}%
\[
\text{ }a_{33}=\frac{d}{dx_{3}}\left[ B_{3}-D_{3}\right] \left( 0\right) , 
\]

\[
b_{11}=\frac{\partial }{\partial x_{1}}\left[ B_{1}-D_{1}\right] \left( \bar{%
x}_{1},0\right) -\frac{\partial h_{1}}{\partial x_{1}}\left( \bar{x}%
_{1},0\right) \text{, }b_{12}=-\frac{\partial D_{1}}{\partial x_{1}}\left( 
\bar{x}_{1},0\right) ,\text{ } 
\]%
\begin{equation}
b_{21}=-\frac{\partial h_{2}}{\partial x_{1}}\left( \bar{x}_{1},0\right) 
\text{, }b_{22}=\frac{d}{dx_{2}}\left[ B_{2}-D_{2}\right] \left( 0\right) -%
\frac{\partial h_{2}}{\partial x_{2}}\left( \bar{x}_{1},0\right) ,\text{ } 
\tag{4.9}
\end{equation}%
\[
b_{31}=\frac{\partial }{\partial x_{1}}\left[ B_{3}-h_{3}\right] \left( \bar{%
x}_{1},0\right) \text{, }b_{33}=\frac{d}{dx_{3}}\left[ B_{3}-D_{3}\right]
\left( \bar{x}_{1},0\right) , 
\]%
\[
c_{11}=\frac{\partial }{\partial x_{1}}\left[ B_{1}-D_{1}\right] \left(
0,0\right) -\frac{\partial h_{1}}{\partial x_{1}}\left( 0,0\right) \text{, }%
c_{12}=-\frac{\partial D_{1}}{\partial x_{1}}\left( 0,0\right) ,\text{ } 
\]%
\[
c_{21}=-\frac{\partial h_{2}}{\partial x_{1}}\left( 0,\bar{x}_{2}\right) 
\text{, }c_{22}=\frac{d}{dx_{2}}\left[ B_{2}-D_{2}\right] \left( 0\right) ,%
\text{ }c_{31}=\frac{\partial }{\partial x_{1}}B_{3}\left( 0,0\right) ,\text{
} 
\]%
\begin{equation}
\text{ }c_{33}=\frac{d}{dx_{3}}\left[ B_{3}-D_{3}\right] \left( 0,0\right) ,
\tag{4.10}
\end{equation}

\[
d_{11}=\frac{\partial }{\partial x_{1}}\left[ B_{1}-D_{1}\right] \left( \bar{%
x}_{1},\bar{x}_{3}\right) -\frac{\partial h_{1}}{\partial x_{1}}\left( \bar{x%
}_{1},\bar{x}_{3}\right) \text{, }d_{12}=-\frac{\partial D_{1}}{\partial
x_{1}}\left( \bar{x}_{1},\bar{x}_{3}\right) ,\text{ }d_{13}= 
\]%
\begin{equation}
-\frac{\partial h_{1}}{\partial x_{3}}\left( \bar{x}_{1},\bar{x}_{3}\right) ,%
\text{ }d_{21}=-\frac{\partial h_{2}}{\partial x_{1}}\left( \bar{x}%
_{1},0\right) \text{, }d_{22}=\frac{d}{dx_{2}}\left[ B_{2}-D_{2}\right]
\left( 0\right) -\frac{\partial h_{2}}{\partial x_{2}}\left( \bar{x}%
_{1},0\right) ,  \tag{4.11}
\end{equation}%
\[
d_{31}=\frac{\partial }{\partial x_{1}}\left[ B_{3}-h_{3}\right] \left( \bar{%
x}_{1},\bar{x}_{3}\right) ,\text{ }d_{33}=\frac{d}{dx_{3}}\left[ B_{3}-D_{3}%
\right] \left( \bar{x}_{1},\bar{x}_{3}\right) , 
\]

\[
k_{11}=\frac{\partial }{\partial x_{1}}\left[ B_{1}-D_{1}\right] \left( \bar{%
x}_{1},0\right) -\frac{\partial h_{1}}{\partial x_{1}}\left( \bar{x}%
_{1},0\right) \text{, }k_{12}=-\frac{\partial D_{1}}{\partial x_{1}}\left( 
\bar{x}_{1},0\right) ,\text{ } 
\]%
\begin{equation}
k_{21}=-\frac{\partial h_{2}}{\partial x_{1}}\left( \bar{x}_{1},\bar{x}%
_{2}\right) \text{, }k_{22}=\frac{d}{dx_{2}}\left[ B_{2}-D_{2}\right] \left( 
\bar{x}_{2}\right) -\frac{\partial h_{2}}{\partial x_{2}}\left( \bar{x}_{1},%
\bar{x}_{2}\right) ,\text{ }  \tag{4.12}
\end{equation}%
\[
k_{31}=\frac{\partial }{\partial x_{1}}B_{3}\left( \bar{x}_{1},0\right) ,%
\text{ }k_{33}=\frac{d}{dx_{3}}\left[ B_{3}-D_{3}\right] \left( \bar{x}%
_{1},0\right) , 
\]

\[
l_{11}=\frac{\partial }{\partial x_{1}}\left[ B_{1}-D_{1}\right] \left( 0,%
\bar{x}_{2}\right) -\frac{\partial h_{1}}{\partial x_{1}}\left( 0,\bar{x}%
_{3}\right) ,\text{ }l_{12}=-\frac{\partial D_{1}}{\partial x_{1}}\left( 0,%
\bar{x}_{2}\right) ,\text{ } 
\]%
\begin{equation}
l_{21}=-\frac{\partial h_{2}}{\partial x_{1}}\left( 0,\bar{x}_{2}\right) ,%
\text{ }l_{22}=\frac{d}{dx_{2}}\left[ B_{2}-D_{2}\right] \left( \bar{x}%
_{2}\right) -\frac{\partial h_{2}}{\partial x_{2}}\left( 0,\bar{x}%
_{2}\right) ,  \tag{4.13}
\end{equation}%
\[
l_{31}=\text{ }\frac{\partial }{\partial x_{1}}\left[ B_{3}\left( 0,\bar{x}%
_{3}\right) -h_{3}\left( 0,\bar{x}_{3}\right) \right] ,\text{ }l_{33}=\frac{%
\partial }{\partial x_{3}}\left[ B_{3}-D_{3}\right] \left( 0,\bar{x}%
_{3}\right) , 
\]%
$\bar{x}_{1}$, $\bar{x}_{2}$ in $\left( 4.9\right) $ and $\left( 4.10\right) 
$\ were defined respectively, by $\ \left( 4.3\right) ,$ $\bar{x}_{1}$, $%
\bar{x}_{3}$ in $\left( 4.11\right) $ were defined by $\left( 4.4\right) ,$ $%
\bar{x}_{1},\bar{x}_{2}$ in $\left( 4.12\right) $\ were defined by $\left(
4.5\right) $ and $\bar{x}_{2},\bar{x}_{3}$ in $\left( 4.13\right) $ were
defined by $\left( 4.6\right) .$

\begin{center}
\textbf{5. local stability analysis of equilibria points }
\end{center}

\bigskip In this section, we derive local stability of the system $\left(
1.1\right) $ at equilibria points $\left( 4.2\right) $. Eigenvalues of\ the
Jacobian matrices $A_{j}$ corresponding to equilibria points $\left(
4.2\right) $ (defined by $\left( 4.7\right) -\left( 4.9\right) $)\ are found
as roots of the equations $\left\vert A_{j}-\lambda \right\vert =0.$

Consider the equilibria point $E_{1}\left( 0,0,0\right) $. Let $a_{ij}$ are
defined by $\left( 4.8\right) .$

\textbf{Theorem 5.1. }Assume\textbf{\ }the assumptions (1)-(5) of Condition
3.1 are satisfied. If $a_{ii}<0$ for $i=1,2,3,$ then the system $\left(
1.1\right) $ is local stabile at the point $E_{1}\left( 0,0,0\right) $; if $%
a_{ii}>0$, then the system $\left( 1.1\right) $ is local unstabile at $%
E_{1}. $

\textbf{Proof. }The\textbf{\ }eigenvalues of\ the Jacobian matrix $A_{1}$\
are found as roots of the equation

\[
\left\vert A_{1}-\lambda \right\vert =\left[ 
\begin{array}{ccc}
a_{11}-\lambda & a_{12} & 0 \\ 
0 & a_{22}-\lambda & 0 \\ 
a_{31} & 0 & a_{33}-\lambda%
\end{array}%
\right] = 
\]%
\[
\left( a_{11}-\lambda \right) \left( a_{22}-\lambda \right) \left(
a_{33}-\lambda \right) =0. 
\]

Hence, $\lambda _{1}=a_{11},$ $\lambda _{2}=a_{22},$ $\lambda _{3}=a_{33}$
are the eigenvalues of the matrix $A_{1}$. By first assumption all
eigenvalues are negative, i.e. the system $\left( 1.1\right) $ is local
stabile at the point $E_{1}$; if $a_{ii}>0$, then the all eigenvalues are
positive, i.e. the system $\left( 1.1\right) $ is local unstabile at $E_{1}.$

Consider the equilibria point $E_{2}\left( \bar{x}_{1},0,0\right) .$ Let $%
b_{ij}$ are defined by $\left( 4.9\right) .$

\bigskip \textbf{Theorem 5.2. }Assume the assumptions (1)-(5) of the
Condition 3.1 are satisfied.\textbf{\ }Let $b_{12}^{2}\leq b_{11}b_{22}.$ If 
$b_{33}<0$ and $b_{11}+b_{22}<0$, then the system $\left( 1.1\right) $ is
local stabile at the point $E_{2}\left( \bar{x}_{1},0,0\right) $; if $%
b_{33}>0$ or $b_{33}\left( b_{11}+b_{22}\right) <0$, then the system $\left(
1.1\right) $ is local unstabile at $E_{2}.$

\bigskip \textbf{Proof. }The eignevalues of\ the Jacobian matrix $A_{2}$\
are found as roots of the equation

\[
\left\vert A_{2}-\lambda \right\vert =\left[ 
\begin{array}{ccc}
b_{11}-\lambda & b_{12} & 0 \\ 
b_{12} & b_{22}-\lambda & 0 \\ 
b_{31} & 0 & b_{33}-\lambda%
\end{array}%
\right] = 
\]%
\[
\left( b_{11}-\lambda \right) \left( b_{22}-\lambda \right) \left(
b_{33}-\lambda \right) -b_{12}^{2}\left( b_{33}-\lambda \right) = 
\]%
\[
\left( b_{33}-\lambda \right) \left[ \left( b_{11}-\lambda \right) \left(
b_{22}-\lambda \right) -b_{12}^{2}\right] =0. 
\]

$\bigskip $Thus, $\lambda _{1}=b_{33},$ $\lambda _{2},$ $\lambda _{3}$ are
the eigenvalues of the matrix $A_{2}$, where $\lambda _{2},$ $\lambda _{3}$
are roots of the equation%
\[
\lambda ^{2}-\left( b_{11}+b_{22}\right) \lambda +b_{11}b_{22}-b_{12}^{2}=0, 
\]%
i.e. \ 
\[
\lambda _{2},\text{ }\lambda _{3}=\frac{\left( b_{11}+b_{22}\right) \pm 
\sqrt{\left( b_{11}+b_{22}\right) ^{2}+4\left(
b_{11}b_{22}-b_{12}^{2}\right) }}{2}. 
\]

That is, if $b_{33}<0$ and $b_{11}+b_{22}<0$, then the all eigenvalues of
the matrix $A_{2}$ are negative, i.e. the system $\left( 1.1\right) $ is
local stabile at the point $E_{2};$ if $b_{33}>0$, $b_{11}+b_{22}>0$ or $%
b_{33}\left( b_{11}+b_{22}\right) <0$, then the all eigenvalues of the
matrix $A_{2}$ are positive, i.e. the system $\left( 1.1\right) $ is local
unstabile at $E_{2}.$

Consider the equilibria point $E_{3}\left( 0,\bar{x}_{2},0\right) .$\ Let $%
c_{ij}$ are defined by $\left( 4.10\right) .$

\textbf{Theorem 5.3. }Assume the assumptions (1)-(5) of the Condition 3.1
are satisfied\textbf{. }Let $c_{12}^{2}\leq c_{11}c_{22},$ $c_{33}<0$ and $%
c_{11}+c_{22}<0$, then the system $\left( 1.1\right) $ is local stabile at
the point $E_{3}\left( 0,\bar{x}_{2},0\right) $; if $c_{33}>0$ or $%
c_{33}\left( c_{11}+c_{22}\right) <0$, then the system $\left( 1.1\right) $
is local unstabile at $E_{3}.$

\bigskip \textbf{Proof. }The\textbf{\ }eigenvalues of\ the Jacobian matrix $%
A_{3}$\ are found as roots

\[
\left\vert A_{3}-\lambda \right\vert =\left[ 
\begin{array}{ccc}
c_{11}-\lambda & c_{12} & 0 \\ 
c_{12} & c_{22}-\lambda & 0 \\ 
c_{31} & 0 & c_{33}-\lambda%
\end{array}%
\right] = 
\]%
\[
\left( c_{11}-\lambda \right) \left( c_{22}-\lambda \right) \left(
c_{33}-\lambda \right) -c_{12}^{2}\left( c_{33}-\lambda \right) = 
\]%
\[
\left( c_{33}-\lambda \right) \left[ \left( c_{11}-\lambda \right) \left(
c_{22}-\lambda \right) -c_{12}^{2}\right] =0. 
\]

$\bigskip $Thus, $\lambda _{1}=c_{33},$ $\lambda _{2},$ $\lambda _{3}$ are
the eigenvalues of the matrix $A_{3}$, where $\lambda _{2},$ $\lambda _{3}$
are roots of the equation%
\[
\lambda ^{2}-\left( c_{11}+c_{22}\right) \lambda +c_{11}c_{22}-c_{12}^{2}=0, 
\]%
i.e. \ 
\[
\lambda _{2},\text{ }\lambda _{3}=\frac{c_{11}+c_{22}\pm \sqrt{\left(
c_{11}+c_{22}\right) ^{2}-4\left( c_{11}c_{22}-c_{12}^{2}\right) }}{2}. 
\]

That is, if $c_{33}<0$ and $c_{11}+c_{22}<0$, then the all eigenvalues of
the matrix $A_{2}$ are negative, i.e. the system $\left( 1.1\right) $ is
local stabile at the point $E_{3};$ if $c_{33}>0$, $c_{11}+c_{22}>0$ or $%
c_{33}\left( c_{11}+c_{22}\right) <0$, then the eigenvalues of the matrix $%
A_{2}$ are positive, i.e. the system $\left( 1.1\right) $ is local unstabile
at $E_{3}.$

Consider the point $E_{4}\left( \bar{x}_{1},0,\bar{x}_{3}\right) .$ Let $%
d_{ij}$ are defined by $\left( 4.11\right) .$

\textbf{Theorem 5.4. }Assume the assumptions (1)-(5) of the Condition 3.1
are satisfied.\textbf{\ }Let $\dsum\limits_{i=1}^{n}d_{ii}<0,$ $%
d_{13}d_{31}d_{22}>-d_{12}^{2}d_{33}$ and $%
d_{11}d_{33}+d_{11}d_{22}+d_{22}d_{33}>d_{12}^{2}+d_{13}d_{31}$. Then the
system $\left( 1.1\right) $ is local stabile at the point $E_{4}\left( \bar{x%
}_{1},0,\bar{x}_{3}\right) $.

\bigskip \textbf{Proof. }Eigenvalues of\ the Jacobian matrix $A_{3}$\ are
found as roots of the equation

\[
\left\vert A_{4}-\lambda \right\vert =\left[ 
\begin{array}{ccc}
d_{11}-\lambda & d_{12} & d_{13} \\ 
d_{12} & d_{22}-\lambda & 0 \\ 
d_{31} & 0 & d_{33}-\lambda%
\end{array}%
\right] = 
\]%
\[
\left( d_{11}-\lambda \right) \left( d_{22}-\lambda \right) \left(
d_{33}-\lambda \right) -d_{12}^{2}\left( d_{33}-\lambda \right)
-d_{13}d_{31}\left( d_{22}-\lambda \right) = 
\]%
\[
\lambda ^{3}-\left( d_{11}+d_{22}+d_{33}\right) \lambda ^{2}+\left(
d_{11}d_{33}+d_{11}d_{22}+d_{22}d_{33}-d_{12}^{2}-d_{13}d_{31}\right)
\lambda + 
\]

\begin{equation}
d_{12}^{2}d_{33}+d_{13}d_{31}d_{22}=0.  \tag{5.1}
\end{equation}%
The roots $\lambda _{1},$ $\lambda _{2}$, $\lambda _{3}$ of $\left(
5.1\right) $ are the eigenvalues of the matrix $A_{4}$. Then by the
fundamental theorem of algebra we have 
\[
\lambda _{1}+\lambda _{2}+\lambda _{3}=d_{11}+d_{22}+d_{33}\text{,} 
\]%
\[
\dsum\limits_{i,j=1}^{3}\lambda _{i}\lambda _{j}=\left(
d_{11}d_{33}+d_{11}d_{22}+d_{22}d_{33}-d_{12}^{2}-d_{13}d_{31}\right) , 
\]%
\[
\lambda _{1}\lambda _{2}\lambda _{3}=-\left[
d_{12}^{2}d_{33}+d_{13}d_{31}d_{22}\right] . 
\]
By the second assumption the all eigenvalues of the matrix $A_{4}$ are
negative, i.e.$\left( 1.1\right) $ is local stabile at $E_{4}\left( \bar{x}%
_{1},0,\bar{x}_{3}\right) .$

Consider the point $E_{5}\left( \bar{x}_{1},\bar{x}_{2},0\right) .$ Let $%
k_{ij}$ are defined by $\left( 4.12\right) .$

\textbf{Theorem 5.5. }Assume the assumptions (1)-(5) of the Condition 3.1
are satisfied.\textbf{\ }Let $k_{12}^{2}\leq k_{11}k_{22}.$ If $k_{33}<0$
and $k_{11}+k_{22}<0$, then the system $\left( 1.1\right) $ is local stabile
at the point $E_{5}\left( \bar{x}_{1},\bar{x}_{2},0\right) $; if $k_{33}>0$
or $k_{33}\left( k_{11}+k_{22}\right) <0$, then the system $\left(
1.1\right) $ is local unstabile at $E_{5}.$

\bigskip \textbf{Proof. }The eigenvalues of\ the Jacobian matrix $A_{5}$\
are found as roots of the equation

\[
\left\vert A_{5}-\lambda \right\vert =\left[ 
\begin{array}{ccc}
k_{11}-\lambda & k_{12} & 0 \\ 
k_{12} & k_{22}-\lambda & 0 \\ 
k_{31} & 0 & k_{33}-\lambda%
\end{array}%
\right] = 
\]%
\[
\left( k_{11}-\lambda \right) \left( k_{22}-\lambda \right) \left(
k_{33}-\lambda \right) -k_{12}^{2}\left( k_{33}-\lambda \right) = 
\]%
\[
\left( k_{33}-\lambda \right) \left[ \left( k_{11}-\lambda \right) \left(
k_{22}-\lambda \right) -k_{12}^{2}\right] =0. 
\]

$\bigskip $Thus, $\lambda _{1}=k_{33},$ $\lambda _{2},$ $\lambda _{3}$ are
the eigenvalues of the matrix $A_{5}$, where $\lambda _{2},$ $\lambda _{3}$
are roots of the equation%
\[
\lambda ^{2}-\left( k_{11}+k_{22}\right) \lambda +k_{11}k_{22}-k_{12}^{2}=0, 
\]%
i.e. \ 
\[
\lambda _{2},\text{ }\lambda _{3}=\frac{k_{11}+k_{22}\pm \sqrt{\left(
k_{11}+k_{22}\right) ^{2}-4\left( k_{11}k_{22}-k_{12}^{2}\right) }}{2}. 
\]

That is, if $k_{33}<0$ and $k_{11}+k_{22}<0$, then the all eigenvalues of
the matrix $A_{2}$ are negative, i.e. the system $\left( 1.1\right) $ is
local stabile at the point $E_{5};$ if $k_{33}>0$, $k_{11}+k_{22}>0$ or $%
k_{33}\left( k_{11}+k_{22}\right) <0$, then the eigenvalues of the matrix $%
A_{2}$ are positive, i.e. the system $\left( 1.1\right) $ is local unstabile
at $E_{5}.$

\bigskip Consider the equilibria point $E_{6}\left( 0,\bar{x}_{2},\bar{x}%
_{3}\right) ,$ where $\bar{x}_{2},$ $\bar{x}_{3}$ is a positive solution of $%
\left( 4.6\right) .$ Let $l_{ij}$ are defined by $\left( 4.13\right) .$

\textbf{Theorem 5.6. }Assume the assumptions (1)-(5) of the Condition 3.1
are satisfied.\textbf{\ }Let $l_{12}l_{21}\leq l_{11}l_{22}.$ If $l_{33}<0$
and $l_{11}+l_{22}<0$, then the system $\left( 1.1\right) $ is local stabile
at the point $E_{6}\left( 0,\bar{x}_{2},\bar{x}_{3}\right) $; if $l_{33}>0$
or $l_{33}\left( l_{11}+l_{22}\right) <0$, then the system $\left(
1.1\right) $ is local unstabile at $E_{6}.$

\bigskip \textbf{Proof. }The eigenvalues of\ the Jacobian matrix $A_{5}$\
are found as roots of the equation

\[
\left\vert A_{6}-\lambda \right\vert =\left[ 
\begin{array}{ccc}
l_{11}-\lambda & l_{12} & 0 \\ 
l_{21} & l_{22}-\lambda & 0 \\ 
l_{31} & 0 & l_{33}-\lambda%
\end{array}%
\right] = 
\]%
\[
\left( l_{11}-\lambda \right) \left( l_{22}-\lambda \right) \left(
l_{33}-\lambda \right) -l_{12}l_{21}\left( l_{33}-\lambda \right) = 
\]%
\[
\left( l_{33}-\lambda \right) \left[ \left( l_{11}-\lambda \right) \left(
l_{22}-\lambda \right) -l_{12}l_{21}\right] =0. 
\]

$\bigskip $Thus, $\lambda _{1}=l_{33},$ $\lambda _{2},$ $\lambda _{3}$ are
the eigenvalues of the matrix $A_{6}$, where $\lambda _{2},$ $\lambda _{3}$
are roots of the equation%
\[
\lambda ^{2}-\left( l_{11}+l_{22}\right) \lambda
+l_{11}l_{22}-l_{12}l_{21}-=0, 
\]

$\bigskip $i.e. \ 
\[
\lambda _{2},\text{ }\lambda _{3}=\frac{l_{11}+l_{22}\pm \sqrt{\left(
l_{11}+l_{22}\right) ^{2}-4\left( l_{11}l_{22}-l_{12}l_{21}\right) }}{2}. 
\]

That is, if $l_{33}<0$ and $l_{11}+l_{22}<0$, then the all eigenvalues of
the matrix $A_{2}$ are negative, i.e. the system $\left( 1.1\right) $ is
local stabile at the point $E_{6};$ if $l_{33}>0$, $l_{11}+l_{22}>0$ or $%
l_{33}\left( l_{11}+l_{22}\right) <0$, then the all eigenvalues of the
matrix $A_{2}$ are positive, i.e. the system $\left( 1.1\right) $ is local
unstabile at $E_{6}.$

\begin{center}
\bigskip \textbf{6. The Lyapunov stability of equilibria points }
\end{center}

\bigskip In this section, we will derive the stability properties of the
system $\left( 1.1\right) $ at points $\left( 4.2\right) $\ in the Lypunov
sense.

Let%
\[
R_{+}^{3}=\left\{ x\in R^{3}\text{: }x_{i}\geq 0,\text{ }i=1,2,3\right\} ,%
\text{ }B_{r}\left( \bar{x}\right) =\left\{ x\in R^{3}\text{, }\left\Vert x-%
\bar{x}\right\Vert _{R^{3}}\leq r^{2}\right\} . 
\]

Let $a_{ij}$ be the real numbers defined by $\left( 4.8\right) $. In this
section we show the following results:

\textbf{Theorem 6.1.} Assume the assumptions (1)-(5) of the Condition 3.1\
are satisfied and $a_{ii}<0$ for $i=1,2,3$. Then the system $\left(
1.1\right) $ is asymptotically stable at the equilibria point $E_{1}\left(
0,0,0\right) $ in the Lyapunov sense.

\textbf{Proof. }Let $A_{1}$ be the linearized matrix with respect to
equilibria point $E_{1}\left( 0,0,0\right) ,$ i.e. 
\[
A_{1}=\left[ 
\begin{array}{ccc}
a_{11} & a_{12} & 0 \\ 
0 & a_{22} & 0 \\ 
a_{31} & 0 & a_{33}%
\end{array}%
\right] \text{, \ }A_{1}^{T}=\left[ 
\begin{array}{ccc}
a_{11} & 0 & a_{31} \\ 
a_{12} & a_{22} & 0 \\ 
0 & 0 & a_{33}%
\end{array}%
\right] . 
\]
We consider the Lyapunov equation 
\[
P_{1}A_{1}+A_{1}^{T}P_{1}=-I,\text{ }P_{1}=\left[ 
\begin{array}{ccc}
p_{11} & p_{12} & p_{13} \\ 
p_{21} & p_{22} & p_{23} \\ 
p_{31} & p_{32} & p_{33}%
\end{array}%
\right] ,\text{ }p_{ij}=p_{ji}, 
\]%
here%
\[
P_{1}A_{1}=\left[ 
\begin{array}{ccc}
p_{11}a_{11}+p_{13}a_{31} & p_{11}a_{12}+p_{12}a_{22} & p_{13}a_{33} \\ 
p_{21}a_{11}+p_{23}a_{31} & p_{21}a_{12}+p_{22}a_{22} & p_{23}a_{33} \\ 
p_{31}a_{11}+p_{33}a_{31} & p_{31}a_{12}+p_{32}a_{22} & p_{33}a_{33}%
\end{array}%
\right] , 
\]

\[
A_{1}^{T}P_{1}=\left[ 
\begin{array}{ccc}
a_{11}p_{11}+a_{31}p_{31} & a_{11}p_{12}+a_{31}p_{32} & 
a_{11}p_{13}+a_{31}p_{33} \\ 
a_{12}p_{11}+a_{22}p_{21} & a_{12}p_{12}+a_{22}p_{22} & 
a_{12}p_{13}+a_{22}p_{23} \\ 
a_{33}p_{31} & a_{33}p_{32} & a_{33}p_{33}%
\end{array}%
\right] , 
\]%
\begin{equation}
P_{1}A_{1}+A_{1}^{T}P_{1}=-I.  \tag{6.1}
\end{equation}

The matrix equation $\left( 6.1\right) $ is equivalent to system of
algebraic equations with respect to $p_{\imath j}$: 
\[
2\left( a_{11}p_{11}+a_{31}p_{13}\right) =-1,\text{ }a_{12}p_{11}+\left(
a_{22}+a_{11}\right) p_{12}+a_{31}p_{23}=0, 
\]%
\[
\left( a_{33}+a_{11}\right) p_{13}+a_{31}p_{33}=0,\text{ }2\left(
a_{12}p_{12}+a_{22}p_{22}\right) =-1, 
\]%
\[
\text{ }\left( a_{22}+a_{33}\right) p_{23}+a_{12}p_{13}=0,\text{ } 
\]%
\[
\text{ }a_{12}p_{13}+\left( a_{22}+a_{33}\right) p_{23}=0,\text{ }%
2p_{33}a_{33}=-1. 
\]

\bigskip\ By solving this system we obtain%
\begin{equation}
\text{ }p_{33}=-\frac{1}{2a_{33}}\text{, }p_{13}=\frac{a_{31}}{2\left(
a_{11}+a_{33}\right) a_{33}},\text{ }p_{11}=-\frac{1}{a_{11}}\left( \frac{1}{%
2}+a_{31}p_{13}\right) ,\text{ }  \tag{6.2}
\end{equation}

\[
p_{23}=-\frac{a_{12}p_{13}}{a_{22}+a_{33}},\text{ }p_{12}=-\frac{\left(
a_{12}p_{11}+a_{31}p_{23}\right) }{\left( a_{11}+a_{22}\right) }\text{, }%
p_{22}=-\frac{-\left( \frac{1}{2}+a_{12}p_{12}\right) }{a_{22}}. 
\]

\ Hence, the eigenvalues of $A_{1}$ are positive if the quadratic function 
\[
V_{1}\left( x\right)
=X^{T}P_{1}X=p_{11}x_{1}^{2}+p_{22}x_{2}^{2}+p_{33}x_{3}^{2}+2p_{12}x_{1}x_{2}+ 
\]%
\[
2p_{13}x_{1}x_{3}+2p_{23}x_{2}x_{3},\text{ }X=\left[ x_{1},x_{2},x_{3}\right]
\]%
is positive defined. It is clear to see that 
\[
V_{1}\left( x\right) =\frac{1}{2}p_{11}x_{1}^{2}+2p_{12}x_{1}x_{2}+\frac{1}{2%
}p_{22}x_{2}^{2}+\frac{1}{2}p_{11}x_{2}^{2}+2p_{13}x_{1}x_{3}+ 
\]%
\begin{equation}
\frac{1}{2}p_{22}x_{2}^{2}+2p_{23}x_{2}x_{3}+p_{33}x_{3}^{2}=  \tag{6.3}
\end{equation}%
\[
\frac{1}{2}p_{11}\left( x_{1}+2\frac{p_{12}}{p_{11}}x_{2}\right) ^{2}+\left( 
\frac{1}{2}p_{22}-2\frac{p_{12}^{2}}{p_{11}}\right) x_{2}^{2}+ 
\]

\[
\frac{1}{2}p_{11}\left( x_{1}+2\frac{p_{12}}{p_{11}}x_{3}\right) ^{2}+\left( 
\frac{1}{2}p_{33}-2\frac{p_{13}^{2}}{p_{11}}\right) x_{3}^{2}+ 
\]%
\[
\frac{1}{2}p_{22}\left( x_{2}+2\frac{p_{23}}{p_{22}}x_{3}\right) ^{2}+\left( 
\frac{1}{2}p_{33}-2\frac{p_{23}^{2}}{p_{22}}\right) x_{3}^{2}>0, 
\]%
when 
\begin{equation}
p_{ii}>0,\text{ }4p_{12}^{2}\leq p_{11}p_{22},\text{ }4p_{13}^{2}\leq
p_{11}p_{33},\text{ }4p_{23}^{2}\leq p_{22}p_{33},  \tag{6.4}
\end{equation}%
\ i.e. the matrix $P_{1}$ is positive defined under the condition $\left(
6.4\right) $. Hence, the quadratic function $V_{1}\left( x\right) $ is a
positive defined Lyapunov function candidate in the certain neighborhood of $%
E_{1}\left( 0,0,0\right) .$ By $\left[ \text{12, Corollary 8.2}\right] $ we
need now to determine a domain $\Omega _{1}$ about the point $E_{1},$ where $%
\dot{V}_{1}\left( x\right) $ is negatively defined and a constant $C$ such
that $\Omega _{C}$ is a subset of $\Omega _{1}$. By assuming $x_{k}\geq 0$, $%
k=1,2,3,$ we will find the solution set of the following inequality 
\begin{equation}
\dot{V}_{1}\left( x\right) =\dsum\limits_{k=1}^{3}\frac{\partial V_{1}}{%
\partial x_{k}}\frac{dx_{k}}{dt}=  \tag{6.5}
\end{equation}%
\[
2\left( p_{11}x_{1}+p_{12}x_{2}+p_{13}x_{3}\right) \left[ B_{1}\left(
x_{1}\right) -D_{1}\left( x_{1},x_{2}\right) -h_{1}\left( x_{1,}x_{3}\right) %
\right] + 
\]%
\[
2\left( p_{12}x_{1}+p_{22}x_{2}+p_{23}x_{3}\right) \left[ B_{2}\left(
x_{2}\right) -D_{2}\left( x_{2}\right) -h_{2}\left( x_{1},x_{2}\right) %
\right] + 
\]%
\[
2\left( p_{13}x_{1}+p_{23}x_{2}+p_{33}x_{3}\right) \left[ B_{3}\left(
x_{1},x_{3}\right) -D_{3}\left( x_{3}\right) -h_{3}\left( x_{1},x_{3}\right) %
\right] \leq 0. 
\]%
Thus, $\left( 6.5\right) $ is satisfied if the following hold%
\[
p_{11}x_{1}+p_{12}x_{2}+p_{13}x_{3}\geq 0,\text{ }%
p_{12}x_{1}+p_{22}x_{2}+p_{23}x_{3}\geq 0,\text{ }%
p_{13}x_{1}+p_{23}x_{2}+p_{33}x_{3}\geq 0, 
\]%
\[
B_{1}\left( x_{1}\right) -D_{1}\left( x_{1},x_{2}\right) -h_{1}\left(
x_{1},x_{3}\right) \leq 0,\text{ }B_{2}\left( x_{2}\right) -D_{2}\left(
x_{2}\right) -h_{2}\left( x_{1},x_{2}\right) \leq 0, 
\]%
\begin{equation}
B_{3}\left( x_{1},x_{3}\right) -D_{3}\left( x_{3}\right) -h_{3}\left(
x_{1},x_{3}\right) \leq 0.  \tag{6.6}
\end{equation}

\textbf{Remark 6.1. }By $\left( 6.2\right) $ the sign of $p_{13}$ is the
same as the sign of $a_{31}$ and the sign of $p_{23}$ is the same as the
sign of $a_{12}a_{31}.$ So, $p_{13}>0$, when $a_{31}>0;$ Hence, $p_{23}>0,$ $%
p_{12}>0$ when $a_{31}>0$ and $a_{12}>0$.\ By assumption $a_{ii}<0$ and $%
\left( 6.2\right) $ we get $p_{11}=-\frac{1}{a_{11}}\left( \frac{1}{2}%
+a_{31}p_{13}\right) >0$, $\ p_{33}>0.$ Since $a_{22}<0$ we get that $%
p_{22}=-\frac{-\left( 1+2a_{12}p_{12}\right) }{2a_{22}}>0$, when $a_{31}>0$
and $\ a_{12}>0$.\ Moreover, by using $\left( 6.2\right) $ we can derive the
conditions on $a_{ij}$ that the assumptions $\left( 6.4\right) $\ are hold.

Here, $b_{ij}$ are real numbers defined by $\left( 4.9\right) .$ Let 
\[
d=\left( b_{11}+b_{33}\right) \left( b_{22}+b_{33}\right) -b_{12}b_{21}, 
\]%
\[
\text{ }D=b_{11}b_{22}\left( b_{11}+b_{22}\right)
-b_{11}b_{12}b_{21}-b_{11}b_{22}b_{12}. 
\]

\bigskip \textbf{Theorem 6.2.} Assume the assumptions (1)-(5) of the
Condition 3.1\ are satisfied. Suppose $b_{ii}<0$ for $i=1,2,3$, $d\neq 0$
and $D\neq 0$. Then the system $\left( 1.1\right) $ is asymptotically stable
at the equilibria point $E_{2}\left( \bar{x}_{1},0,0\right) $ in the
Lyapunov sense.

\textbf{Proof. }Let $A_{2}$ be the linearized matrix with respect to
equilibria point $E_{2}\left( \bar{x}_{1},0,0\right) ,$ i.e. 
\[
A_{2}=\left[ 
\begin{array}{ccc}
b_{11} & b_{12} & 0 \\ 
b_{21} & b_{22} & 0 \\ 
b_{31} & 0 & b_{33}%
\end{array}%
\right] ,\text{ }A_{2}^{T}=\left[ 
\begin{array}{ccc}
b_{11} & b_{21} & b_{31} \\ 
b_{12} & b_{22} & 0 \\ 
0 & 0 & b_{33}%
\end{array}%
\right] . 
\]
We consider the Lyapunov equation 
\begin{equation}
P_{2}A_{2}+A_{2}^{T}P_{2}=-I,\text{ }P_{2}=\left[ 
\begin{array}{ccc}
p_{11} & p_{12} & p_{13} \\ 
p_{21} & p_{22} & p_{23} \\ 
p_{31} & p_{32} & p_{33}%
\end{array}%
\right] ,\text{ }p_{ij}=p_{ji},  \tag{6.7}
\end{equation}%
where%
\[
P_{2}A_{2}=\left[ 
\begin{array}{ccc}
p_{11}b_{11}+p_{12}b_{21}+p_{13}b_{31} & p_{11}b_{12}+p_{12}b_{22} & 
p_{13}b_{33} \\ 
p_{21}b_{11}+p_{22}b_{21}+p_{23}b_{31} & p_{21}b_{12}+p_{22}b_{22} & 
p_{23}b_{33} \\ 
p_{31}b_{11}+p_{32}b_{21}+p_{33}b_{31} & p_{31}b_{12}+p_{32}b_{22} & 
p_{33}b_{33}%
\end{array}%
\right] 
\]

\[
A_{2}^{T}P_{2}=\left[ 
\begin{array}{ccc}
b_{11}p_{11}+b_{21}p_{21}+b_{31}p_{31} & 
b_{11}p_{12}+b_{21}p_{22}+b_{31}p_{32} & 
b_{11}p_{13}+b_{21}p_{23}+b_{31}p_{33} \\ 
b_{12}p_{11}+b_{22}p_{21} & b_{12}p_{12}+b_{22}p_{22} & 
b_{12}p_{13}+b_{22}p_{23} \\ 
b_{33}p_{31} & b_{33}p_{32} & b_{33}p_{33}%
\end{array}%
\right] , 
\]%
\begin{equation}
P_{2}A_{2}+A_{2}^{T}P_{2}=-I.  \tag{6.8}
\end{equation}

The matrix equation $\left( 6.1\right) $ is equivalent to system of
algebraic equations with respect to $p_{\imath j}$ 
\[
2\left( b_{11}p_{11}+b_{21}p_{12}+b_{31}p_{13}\right) =-1,\text{ }%
b_{12}p_{11}+\left( b_{22}+b_{11}\right) p_{12}+b_{21}p_{22}+ 
\]%
\[
b_{31}p_{23}=0,\text{\ }\left( b_{33}+b_{11}\right)
p_{13}+b_{21}p_{23}+b_{31}p_{33}=0,\text{\ } 
\]%
\[
\text{ }2\left( b_{12}p_{12}+p_{22}b_{22}\right) =-1,\text{ }\left(
b_{33}+b_{22}\right) p_{23}+b_{12}p_{13}=0, 
\]%
\[
\text{ }\left( b_{11}+b_{33}\right) p_{13}+b_{21}p_{23}+b_{31}p_{33}=0, 
\]%
\[
b_{12}p_{13}+\left( b_{22}+b_{33}\right) p_{23}=0,\text{ }2p_{33}b_{33}=-1. 
\]

By solving this system we obtain 
\[
p_{33}=-\frac{1}{2b_{33}},\text{ }p_{13}=\frac{d_{1}}{d},\text{ }p_{23}=%
\frac{d_{2}}{d},\text{ }p_{11}=\frac{D_{1}}{D},\text{ }p_{12}=\frac{D_{2}}{D}%
,\text{ }p_{22}=\frac{D_{3}}{D}, 
\]%
where 
\begin{equation}
\text{ }d_{1}=-\frac{b_{21}b_{31}}{2b_{33}},\text{ }d_{2}=\frac{b_{31}}{%
2b_{33}}\left( b_{11}+b_{33}\right) ,  \tag{6.9}
\end{equation}%
\[
\text{ }D_{1}=-\frac{1}{2}b_{21}^{2}+b_{22}\left( b_{11}+b_{22}\right)
\left( \frac{1}{2}+b_{31}p_{13}\right) + 
\]%
\[
\left( \frac{1}{2}+b_{31}p_{13}\right)
b_{12}b_{21}+b_{21}b_{22}b_{31}p_{23}, 
\]%
\[
\text{ }D_{2}=\frac{1}{2}b_{11}b_{21}+b_{12}b_{22}\left( \frac{1}{2}%
+b_{31}p_{13}\right) -b_{11}b_{22}b_{31}p_{23}, 
\]%
\[
\text{ }D_{3}=b_{11}b_{12}b_{31}p_{23}+\frac{1}{2}b_{12}b_{21}-\frac{1}{2}%
b_{11}\left( b_{11}+b_{22}\right) -b_{12}^{2}\left( \frac{1}{2}%
+b_{31}p_{13}\right) . 
\]%
Hence, the eigenvalues of $A_{2}$ are positive if the quadratic function 
\[
V_{2}\left( x\right)
=X^{T}P_{2}X=p_{11}x_{1}^{2}+p_{22}x_{2}^{2}+p_{33}x_{3}^{2}+2p_{12}x_{1}x_{2}+ 
\]%
\[
2p_{13}x_{1}x_{3}+2p_{23}x_{2}x_{3} 
\]%
is positive defined. By assumption we get that $p_{33}>0.$ Moreover, $%
p_{kk}>0$ for $k=1,2$, when $\frac{D_{1}}{D}>0,$ $\frac{D_{3}}{D}>0.$ Hence,
in a similar way we obtain that $V_{2}\left( x\right) $\ is positive
defined, if $\frac{D_{1}}{D}>0,$ $\frac{D_{3}}{D}>0$ and when the estimate
of type $\left( 6.4\right) $ is satisfied.

By reasoning as in the proof of Theorem 6.1 we obtain that the inequality 
\begin{equation}
\dot{V}_{2}\left( x\right) =\dsum\limits_{k=1}^{3}\frac{\partial V_{2}}{%
\partial x_{k}}\frac{dx_{k}}{dt}\leq 0  \tag{6.10}
\end{equation}%
is valid if the following holds%
\[
p_{11}x_{1}+p_{12}x_{2}+p_{13}x_{3}\geq 0,\text{ }%
p_{12}x_{1}+p_{22}x_{2}+p_{23}x_{3}\geq 0,\text{ }%
p_{13}x_{1}+p_{23}x_{2}+p_{33}x_{3}\geq 0, 
\]%
\[
B_{1}\left( x_{1}\right) -D_{1}\left( x_{1},x_{2}\right) -h_{1}\left(
x_{1},x_{3}\right) \leq 0,\text{ }B_{2}\left( x_{2}\right) -D_{2}\left(
x_{2}\right) -h_{2}\left( x_{1},x_{2}\right) \leq 0, 
\]%
\begin{equation}
B_{3}\left( x_{1},x_{3}\right) -D_{3}\left( x_{3}\right) -h_{3}\left(
x_{1},x_{3}\right) \leq 0.  \tag{6.11}
\end{equation}

\textbf{Remark 6.2. }In view of\textbf{\ }$\left( 6.2\right) $, $p_{kk}>0$
when $\left( b_{12}^{2}-b_{11}b_{22}\right) <0,$ $D_{1}>0,$ $D_{3}>0$ or $%
\left( b_{12}^{2}-b_{11}b_{22}\right) >0,$ $D_{1}<0,$ $D_{3}<0.$ Moreover,
by using $\left( 6.9\right) $ we can derived the conditions on $b_{ij}$ that
the assumptions of type $\left( 6.4\right) $\ are hold.

Here, $c_{ij}$ are real numbers\ defined by $\left( 4.10\right) .$ Let%
\[
d=\left( c_{11}+c_{33}\right) \left( c_{22}+c_{33}\right) -c_{12}c_{21},%
\text{ } 
\]%
\[
D=c_{11}c_{22}\left( c_{11}+c_{22}\right)
-c_{11}c_{12}c_{21}-c_{11}c_{22}c_{12}. 
\]

\textbf{Theorem 6.3.} Assume the assumptions (1)-(5) of the Condition 3.1\
are satisfied. Suppose $c_{ii}<0$ for $i=1,2,3$, $d\neq 0$ and $D\neq 0$.
Then the system $\left( 1.1\right) $ is asymptotically stable at the
equilibria point $E_{2}\left( \bar{x}_{1},0,0\right) $ in the Lyapunov sense.

\textbf{Proof. }Let $A_{3}$ be the linearized matrix with respect to
equilibria point $E_{3}\left( 0,\bar{x}_{2},0\right) ,$ i.e. 
\[
A_{3}=\left[ 
\begin{array}{ccc}
c_{11} & c_{12} & 0 \\ 
c_{12} & c_{22} & 0 \\ 
c_{31} & 0 & c_{33}%
\end{array}%
\right] ,\text{ }A_{3}^{T}=\left[ 
\begin{array}{ccc}
c_{11} & c_{12} & c_{31} \\ 
c_{12} & c_{22} & 0 \\ 
0 & 0 & c_{33}%
\end{array}%
\right] . 
\]
We consider the Lyapunov equation 
\begin{equation}
P_{3}A_{3}+A_{3}^{T}P_{3}=-I,\text{ }P_{3}=\left[ 
\begin{array}{ccc}
p_{11} & p_{12} & p_{13} \\ 
p_{21} & p_{22} & p_{23} \\ 
p_{31} & p_{32} & p_{33}%
\end{array}%
\right] ,\text{ }p_{ij}=p_{ji}.  \tag{6.12}
\end{equation}

\bigskip By solving $\left( 6.12\right) $ as in the Theorem 6.2 we obtain 
\[
p_{33}=-\frac{1}{2c_{33}},\text{ }p_{13}=\frac{d_{1}}{d},\text{ }p_{23}=%
\frac{d_{2}}{d},\text{ }p_{11}=\frac{D_{1}}{D},\text{ }p_{12}=\frac{D_{2}}{D}%
,\text{ }p_{22}=\frac{D_{3}}{D}, 
\]%
where%
\begin{equation}
\text{ }d_{1}=-\frac{c_{21}c_{31}}{2c_{33}},\text{ }d_{2}=\frac{c_{31}}{%
2c_{33}}\left( c_{11}+c_{33}\right) ,  \tag{6.13}
\end{equation}%
\[
\text{ }D_{1}=-\frac{1}{2}c_{21}^{2}+c_{22}\left( c_{11}+c_{22}\right)
\left( \frac{1}{2}+c_{31}p_{13}\right) + 
\]%
\[
\left( \frac{1}{2}+c_{31}p_{13}\right)
c_{12}c_{21}+c_{21}c_{22}c_{31}p_{23}, 
\]%
\[
\text{ }D_{2}=\frac{1}{2}c_{11}c_{21}+c_{12}c_{22}\left( \frac{1}{2}%
+c_{31}p_{13}\right) -c_{11}c_{22}c_{31}p_{23}, 
\]%
\[
\text{ }D_{3}=c_{11}c_{12}c_{31}p_{23}+\frac{1}{2}c_{12}c_{21}-\frac{1}{2}%
c_{11}\left( c_{11}+c_{22}\right) -c_{12}^{2}\left( \frac{1}{2}%
+c_{31}p_{13}\right) . 
\]

\bigskip Hence, the eigenvalues of $A_{3}$ are positive if the quadratic
function 
\[
V_{3}\left( x\right)
=X^{T}P_{2}X=p_{11}x_{1}^{2}+p_{22}x_{2}^{2}+p_{33}x_{3}^{2}+2p_{12}x_{1}x_{2}+ 
\]%
\[
2p_{13}x_{1}x_{3}+2p_{23}x_{2}x_{3} 
\]%
is positive defined. In a similar way we obtain that $V_{3}\left( x\right) $%
\ is positive defined, when $\frac{D_{1}}{D}>0,$ $\frac{D_{3}}{D}>0$ and the
conditions of type $\left( 6.4\right) $ are hold.

By reasoning as in the proof of Theorem 6.1 we obtain that the inequality 
\[
\dot{V}_{3}\left( x\right) =\dsum\limits_{k=1}^{3}\frac{\partial V_{3}}{%
\partial x_{k}}\frac{dx_{k}}{dt}\leq 0 
\]%
is valid if the following are hold%
\[
p_{11}x_{1}+p_{12}x_{2}+p_{13}x_{3}\geq 0,\text{ }%
p_{12}x_{1}+p_{22}x_{2}+p_{23}x_{3}\geq 0,\text{ }%
p_{13}x_{1}+p_{23}x_{2}+p_{33}x_{3}\geq 0, 
\]%
\[
B_{1}\left( x_{1}\right) -D_{1}\left( x_{1},x_{2}\right) -h_{1}\left(
x_{1},x_{3}\right) \leq 0,\text{ }B_{2}\left( x_{2}\right) -D_{2}\left(
x_{2}\right) -h_{2}\left( x_{1},x_{2}\right) \leq 0, 
\]%
\begin{equation}
B_{3}\left( x_{1},x_{3}\right) -D_{3}\left( x_{3}\right) -h_{3}\left(
x_{1},x_{3}\right) \leq 0.  \tag{6.14}
\end{equation}

\bigskip \textbf{Remark 6.3. }By $\left( 6.13\right) $, $p_{kk}>0$ when $%
\left( c_{12}^{2}-c_{11}c_{22}\right) <0,$ $D_{1}>0,$ $D_{3}>0$ or $\left(
c_{12}^{2}-c_{11}c_{22}\right) >0,$ $D_{1}<0,$ $D_{3}<0.$ Moreover, by using 
$\left( 6.13\right) $ we can derived the conditions on $c_{ij}$ that the
assumptions of type $\left( 6.4\right) $\ are hold.

Consider the stable point $E_{4}\left( \bar{x}_{1},0,\bar{x}_{3}\right) .$
Here, $d_{ij}$ are real numbers\ defined by $\left( 4.11\right) .$\ Let 
\[
d=\left( d_{11}+d_{33}\right) \left( d_{22}+d_{33}\right) -d_{12}d_{21}, 
\]%
\[
\text{ }D=d_{11}d_{22}\left( d_{11}+d_{22}\right)
-d_{11}d_{12}d_{21}-d_{11}d_{22}d_{12}. 
\]

\textbf{Theorem 6.4.} Assume the assumptions (1)-(5) of the Condition 3.1\
are satisfied. Suppose $d_{ii}<0$ for $i=1,2,3$, $d\neq 0$ and $D\neq 0$.
Then the system $\left( 1.1\right) $ is asymptotically stable at the
equilibria point $E_{4}\left( \bar{x}_{1},0,\bar{x}_{3}\right) $ in the
Lyapunov sense.

\textbf{Proof. }Let $A_{4}$ be the linearized matrix with respect to
equilibria point $E_{4}\left( \bar{x}_{1},0,\bar{x}_{3}\right) ,$ i.e.

\[
\text{ }A_{4}=\left[ 
\begin{array}{ccc}
d_{11} & d_{12} & d_{13} \\ 
d_{21} & d_{22} & 0 \\ 
d_{31} & 0 & d_{33}%
\end{array}%
\right] ,\text{ }A_{4}^{T}=\left[ 
\begin{array}{ccc}
d_{11} & d_{21} & d_{31} \\ 
d_{12} & d_{22} & 0 \\ 
d_{13} & 0 & d_{33}%
\end{array}%
\right] . 
\]
We consider the Lyapunov equation 
\begin{equation}
P_{4}A_{4}+A_{4}^{T}P_{4}=-I,\text{ }P_{4}=\left[ 
\begin{array}{ccc}
p_{11} & p_{12} & p_{13} \\ 
p_{21} & p_{22} & p_{23} \\ 
p_{31} & p_{32} & p_{33}%
\end{array}%
\right] ,\text{ }p_{ij}=p_{ji}.  \tag{6.15}
\end{equation}%
It is clear that%
\[
P_{4}A_{4}=\left[ 
\begin{array}{ccc}
d_{11}p_{11}+d_{21}p_{12}+d_{31}p_{13} & d_{12}p_{11}+d_{22}p_{12} & 
d_{13}p_{11}+d_{33}p_{13} \\ 
d_{11}p_{21}+d_{21}p_{22}+d_{31}p_{23} & d_{12}p_{21}+d_{22}p_{22} & 
d_{13}p_{21}+d_{33}p_{23} \\ 
d_{11}p_{31}+d_{21}p_{32}+d_{31}p_{33} & d_{12}p_{31}+d_{22}p_{32} & 
d_{13}p_{31}+d_{33}p_{33}%
\end{array}%
\right] , 
\]%
\[
A_{4}^{T}P_{4}=\left[ 
\begin{array}{ccc}
d_{11}p_{11}+d_{21}p_{21}+d_{31}p_{31} & 
d_{11}p_{12}+d_{21}p_{22}+d_{31}p_{32} & 
d_{11}p_{13}+d_{21}p_{23}+d_{31}p_{33} \\ 
d_{12}p_{11}+d_{22}p_{21} & d_{12}p_{12}+d_{22}p_{22} & 
d_{12}p_{13}+d_{22}p_{23} \\ 
d_{13}p_{11}+d_{33}p_{31} & d_{13}p_{12}+d_{33}p_{32} & 
d_{13}p_{13}+d_{33}p_{33}%
\end{array}%
\right] , 
\]

\[
P_{4}A_{4}+A_{4}^{T}P_{4}=\left[ 
\begin{array}{ccc}
d_{11}p_{11}+d_{21}p_{12}+d_{31}p_{13}+d_{11}p_{11}+d_{21}p_{21}+d_{31}p_{31}
& d_{21}p_{11}+d_{22}p_{12}+d_{11}p_{12}+d_{21}p_{22}+d_{31}p_{32} & 
d_{13}p_{11}+d_{33}p_{13}+d_{11}p_{13}+d_{21}p_{23}+d_{31}p_{33} \\ 
d_{11}p_{21}+d_{21}p_{22}+d_{31}p_{23}+d_{12}p_{11}+d_{22}p_{21} & 
d_{12}p_{21}+d_{22}p_{22}+d_{12}p_{12}+d_{22}p_{22} & 
d_{13}p_{21}+d_{33}p_{23}+d_{12}p_{13}+d_{22}p_{23} \\ 
d_{11}p_{31}+d_{21}p_{32}+d_{31}p_{33}+d_{13}p_{11}+d_{33}p_{31} & 
d_{12}p_{31}+d_{22}p_{32}+d_{13}p_{12}+d_{33}p_{32} & 
d_{13}p_{31}+d_{33}p_{33}+d_{13}p_{13}+d_{33}p_{33}%
\end{array}%
\right] . 
\]

From $\left( 6.15\right) $ we obtain the following system of the equations
in $p_{ij}:$

\[
2\left( d_{11}p_{11}+d_{21}p_{12}+d_{31}p_{13}\right) =-1,\text{ }%
d_{21}p_{11}+\left( d_{22}+d_{11}\right) p_{12}+d_{21}p_{22}+d_{31}p_{23}=0, 
\]%
\[
\text{ }d_{13}p_{11}+\left( d_{33}+d_{11}\right)
p_{13}+d_{21}p_{23}+d_{31}p_{33}=0,\text{ }2\left(
d_{12}p_{12}+d_{22}p_{22}\right) =-1, 
\]%
\[
d_{12}p_{13}+\left( d_{33}+d_{22}\right) p_{23}+d_{13}p_{12}=0,\text{ }%
2\left( d_{13}p_{13}+d_{33}p_{33}\right) =-1. 
\]

By taking%
\[
p_{22}=-\frac{1}{d_{22}}\left( \frac{1}{2}+d_{12}p_{12}\right) ,\text{ }%
p_{33}=-\frac{1}{d_{33}}\left( \frac{1}{2}+d_{13}p_{13}\right) 
\]%
\ in \ the other equations we get 
\[
2\left( d_{11}p_{11}+d_{21}p_{12}+d_{31}p_{13}\right) =-1,\text{ } 
\]%
\begin{equation}
d_{21}p_{11}+\left( d_{22}+d_{11}-\frac{d_{12}d_{21}}{d_{22}}\right)
p_{12}+d_{31}p_{23}=\frac{d_{12}}{2d_{22}},  \tag{6.16}
\end{equation}%
\[
d_{13}p_{11}+\left( d_{33}+d_{11}-\frac{d_{13}d_{31}}{d_{33}}\right)
p_{13}+d_{21}p_{23}=\frac{d_{13}}{2d_{33}} 
\]%
\[
d_{12}p_{13}+\left( d_{33}+d_{22}\right) p_{23}+d_{13}p_{12}=0. 
\]%
By solving the system $\left( 6.16\right) $ we get 
\[
p_{11}=\frac{D_{1}}{D},\text{ }p_{12}=\frac{D_{2}}{D},\text{ }p_{13}=\frac{%
D_{3}}{D},\text{ }p_{23}=\frac{D_{4}}{D}, 
\]%
where

\[
D=\left\vert 
\begin{array}{cccc}
2d_{11} & 2d_{21} & 2d_{31} & 0 \\ 
d_{21} & d_{0} & 0 & d_{31} \\ 
0 & d_{13} & d_{12} & d_{22}+d_{33} \\ 
0 & d_{13} & d_{12} & d_{22}+d_{33}%
\end{array}%
\right\vert , 
\]

\[
D_{1}=\left\vert 
\begin{array}{cccc}
-1 & 2d_{21} & 2d_{31} & 0 \\ 
\frac{d_{12}}{2d_{22}} & d_{0} & 0 & d_{31} \\ 
\frac{d_{13}}{2d_{33}} & d_{13} & d_{12} & d_{22}+d_{33} \\ 
0 & d_{13} & d_{12} & d_{22}+d_{33}%
\end{array}%
\right\vert ,\text{ }D_{2}=\left\vert 
\begin{array}{cccc}
2d_{11} & -1 & 2d_{31} & 0 \\ 
d_{21} & \frac{d_{12}}{2d_{22}} & 0 & d_{31} \\ 
0 & \frac{d_{13}}{2d_{33}} & d_{12} & d_{22}+d_{33} \\ 
0 & 0 & d_{12} & d_{22}+d_{33}%
\end{array}%
\right\vert ,\text{ } 
\]%
\[
D_{3}=%
\begin{array}{cccc}
2d_{11} & 2d_{21} & -1 & 0 \\ 
d_{21} & d_{0} & \frac{d_{12}}{2d_{22}} & d_{31} \\ 
0 & d_{13} & \frac{d_{13}}{2d_{33}} & d_{22}+d_{33} \\ 
0 & d_{13} & 0 & d_{22}+d_{33}%
\end{array}%
\text{ , \ }D_{4}=\left\vert 
\begin{array}{cccc}
2d_{11} & 2d_{21} & 2d_{31} & -1 \\ 
d_{21} & d_{0} & 0 & \frac{d_{12}}{2d_{22}} \\ 
0 & d_{13} & d_{12} & \frac{d_{13}}{2d_{33}} \\ 
0 & d_{13} & d_{12} & 0%
\end{array}%
\right\vert ; 
\]%
here, 
\begin{equation}
d_{0}=d_{22}+d_{11}-\frac{d_{12}d_{21}}{d_{22}},\text{ }b_{0}=d_{33}+d_{11}-%
\frac{d_{13}d_{31}}{d_{33}},  \tag{6.17}
\end{equation}%
\[
p_{22}=-\frac{1}{d_{22}}\left( \frac{1}{2}+d_{12}p_{12}\right) =-\frac{1}{%
d_{22}}\left( \frac{1}{2}+d_{12}\frac{D_{2}}{D}\right) ,\text{ }p_{33}=-%
\frac{1}{d_{33}}\left( \frac{1}{2}+d_{13}\frac{D_{3}}{D}\right) . 
\]

Thus, the eigenvalues of $A_{4}$ are positive if the quadratic function 
\[
V_{4}\left( x\right)
=X^{T}P_{2}X=p_{11}x_{1}^{2}+p_{22}x_{2}^{2}+p_{33}x_{3}^{2}+2p_{12}x_{1}x_{2}+ 
\]%
\[
2p_{13}x_{1}x_{3}+2p_{23}x_{2}x_{3} 
\]%
is positive defined. In a similar way we obtain that $V_{4}\left( x\right) $%
\ is positive defined, when the conditions of type $\left( 6.4\right) $ are
hold.

By reasoning as in the proof of Theorem 6.1 we obtain that the inequality 
\[
\dot{V}_{4}\left( x\right) =\dsum\limits_{k=1}^{3}\frac{\partial V_{4}}{%
\partial x_{k}}\frac{dx_{k}}{dt}\leq 0 
\]%
is valid if the following are satisfied%
\[
p_{11}x_{1}+p_{12}x_{2}+p_{13}x_{3}\geq 0,\text{ }%
p_{12}x_{1}+p_{22}x_{2}+p_{23}x_{3}\geq 0,\text{ }%
p_{13}x_{1}+p_{23}x_{2}+p_{33}x_{3}\geq 0, 
\]%
\[
B_{1}\left( x_{1}\right) -D_{1}\left( x_{1},x_{2}\right) -h_{1}\left(
x_{1},x_{3}\right) \leq 0,\text{ }B_{2}\left( x_{2}\right) -D_{2}\left(
x_{2}\right) -h_{2}\left( x_{1},x_{2}\right) \leq 0, 
\]%
\begin{equation}
B_{3}\left( x_{1},x_{3}\right) -D_{3}\left( x_{3}\right) -h_{3}\left(
x_{1},x_{3}\right) \leq 0.  \tag{6.18}
\end{equation}

\bigskip \textbf{Remark 6.4. \ }By $\left( 6.17\right) $, $p_{kk}>0$ when $%
\frac{D_{1}}{D}>0,$ $-\frac{1}{d_{22}}\left( \frac{1}{2}+d_{12}\frac{D_{2}}{D%
}\right) >0,$ $-\frac{1}{d_{33}}\left( \frac{1}{2}+d_{13}\frac{D_{3}}{D}%
\right) >0.$ Moreover, by using $\left( 6.17\right) $ we can derived the
conditions on $d_{ij}$ that the assumptions of type $\left( 6.4\right) $\
are hold.

Here, $k_{ij}$ are real numbers defined by $\left( 4.12\right) .$ Let%
\[
d=\left( k_{11}+k_{33}\right) \left( k_{22}+k_{33}\right) -k_{12}k_{21},%
\text{ } 
\]%
\[
D=k_{11}k_{22}\left( k_{11}+k_{22}\right)
-k_{11}k_{12}k_{21}-k_{11}k_{22}k_{12}. 
\]

\textbf{Theorem 6.5.} Assume the assumptions (1)-(5) of the Condition 3.1\
are satisfied. Suppose $k_{ii}<0$ for $i=1,2,3$, $d\neq 0$ and $D\neq 0$.
Then the system $\left( 1.1\right) $ is asymptotically stable at the
equilibria point $E_{5}\left( \bar{x}_{1},\bar{x}_{2},0\right) $ in the
Lyapunov sense.

\textbf{Proof. }Let $A_{5}$ be the linearized matrix with respect to
equilibria point $E_{5}\left( \bar{x}_{1},\bar{x}_{2},0\right) ,$ i.e. 
\[
A_{5}=\left[ 
\begin{array}{ccc}
k_{11} & k_{12} & 0 \\ 
k_{21} & k_{22} & 0 \\ 
k_{31} & 0 & k_{33}%
\end{array}%
\right] ,\text{ }A_{5}^{T}=\left[ 
\begin{array}{ccc}
k_{11} & k_{21} & k_{31} \\ 
k_{12} & k_{22} & 0 \\ 
0 & 0 & k_{33}%
\end{array}%
\right] . 
\]%
We consider the Lyapunov equation 
\begin{equation}
P_{5}A_{5}+A_{5}^{T}P_{5}=-I,\text{ }P_{5}=\left[ 
\begin{array}{ccc}
p_{11} & p_{12} & p_{13} \\ 
p_{21} & p_{22} & p_{23} \\ 
p_{31} & p_{32} & p_{33}%
\end{array}%
\right] ,\text{ }p_{ij}=p_{ji}.  \tag{6.19}
\end{equation}

\bigskip By solving $\left( 6.19\right) ,$ in a similar way as in the
Theorem 6.2 we obtain 
\[
p_{33}=-\frac{1}{2k_{33}},\text{ }p_{13}=\frac{d_{1}}{d},\text{ }p_{23}=%
\frac{d_{2}}{d},\text{ }p_{11}=\frac{D_{1}}{D},\text{ }p_{12}=\frac{D_{2}}{D}%
,\text{ }p_{22}=\frac{D_{3}}{D}, 
\]%
where 
\begin{equation}
\text{ }d_{1}=-\frac{k_{21}k_{31}}{2k_{33}},\text{ }d_{2}=\frac{k_{31}}{%
2k_{33}}\left( k_{11}+k_{33}\right) ,  \tag{6.20}
\end{equation}%
\[
\text{ }D_{1}=-\frac{1}{2}k_{21}^{2}+k_{22}\left( k_{11}+k_{22}\right)
\left( \frac{1}{2}+k_{31}p_{13}\right) + 
\]%
\[
\left( \frac{1}{2}+k_{31}p_{13}\right)
k_{12}k_{21}+k_{21}k_{22}k_{31}p_{23}, 
\]%
\[
\text{ }D_{2}=\frac{1}{2}k_{11}k_{21}+k_{12}k_{22}\left( \frac{1}{2}%
+k_{31}p_{13}\right) -k_{11}k_{22}k_{31}p_{23}, 
\]%
\[
\text{ }D_{3}=k_{11}k_{12}k_{31}p_{23}+\frac{1}{2}k_{12}k_{21}-\frac{1}{2}%
k_{11}\left( k_{11}+k_{22}\right) -k_{12}^{2}\left( \frac{1}{2}%
+k_{31}p_{13}\right) . 
\]

\bigskip Hence, the eigenvalues of $A_{5}$ are positive if the quadratic
function 
\[
V_{5}\left( x\right)
=X^{T}P_{2}X=p_{11}x_{1}^{2}+p_{22}x_{2}^{2}+p_{33}x_{3}^{2}+2p_{12}x_{1}x_{2}+ 
\]%
\[
2p_{13}x_{1}x_{3}+2p_{23}x_{2}x_{3} 
\]%
is positive defined. In a similar way we obtain that $V_{5}\left( x\right) $%
\ is positive defined, when $\frac{D_{1}}{D}>0,$ $\frac{D_{3}}{D}>0$ and the
conditions of type $\left( 6.4\right) $ are satisfied. By reasoning as in
the proof of Theorem 6.1 we obtain that the inequality 
\[
\dot{V}_{5}\left( x\right) =\dsum\limits_{k=1}^{3}\frac{\partial V_{5}}{%
\partial x_{k}}\frac{dx_{k}}{dt}\leq 0 
\]%
is valid if the following holds%
\[
p_{11}x_{1}+p_{12}x_{2}+p_{13}x_{3}\geq 0,\text{ }%
p_{12}x_{1}+p_{22}x_{2}+p_{23}x_{3}\geq 0,\text{ }%
p_{13}x_{1}+p_{23}x_{2}+p_{33}x_{3}\geq 0, 
\]%
\[
B_{1}\left( x_{1}\right) -D_{1}\left( x_{1},x_{2}\right) -h_{1}\left(
x_{1},x_{3}\right) \leq 0,\text{ }B_{2}\left( x_{2}\right) -D_{2}\left(
x_{2}\right) -h_{2}\left( x_{1},x_{2}\right) \leq 0, 
\]%
\begin{equation}
B_{3}\left( x_{1},x_{3}\right) -D_{3}\left( x_{3}\right) -h_{3}\left(
x_{1},x_{3}\right) \leq 0.  \tag{6.21}
\end{equation}

\bigskip \textbf{Remark 6.5. }In view of $\left( 6.17\right) $, $p_{kk}>0$
when $\left( k_{12}^{2}-k_{11}k_{22}\right) <0,$ $D_{1}>0,$ $D_{3}>0$ or $%
\left( k_{12}^{2}-k_{11}k_{22}\right) >0,$ $D_{1}<0,$ $D_{3}<0.$ Moreover,
by using $\left( 6.20\right) $ we can derived the conditions on $k_{ij}$
that the assumptions of type $\left( 6.4\right) $\ are hold.

Here, $l_{ij}$ are real numbers\ defined by $\left( 4.13\right) .$ Let%
\[
d=\left( l_{11}+l_{33}\right) \left( l_{22}+l_{33}\right) -l_{12}l_{21},%
\text{ } 
\]%
\[
D=l_{11}l_{22}\left( l_{11}+l_{22}\right)
-l_{11}l_{12}l_{21}-l_{11}l_{22}l_{12}. 
\]

\textbf{Theorem 6.6.} Assume the assumptions (1)-(5) of the Condition 3.1\
are satisfied. Suppose $l_{ii}<0$ for $i=1,2,3$, $d\neq 0$ and $D\neq 0$.
Then the system $\left( 1.1\right) $ is asymptotically stable at the
equilibria point $E_{6}\left( 0,\bar{x}_{2},\bar{x}_{3}\right) $ in the
Lyapunov sense.

\textbf{Proof. }Let $A_{6}$ be the linearized matrix with respect to
equilibria point $E_{6}\left( 0,\bar{x}_{2},\bar{x}_{3}\right) ,$ i.e. 
\[
A_{6}=\left[ 
\begin{array}{ccc}
l_{11} & l_{12} & 0 \\ 
l_{21} & l_{22} & 0 \\ 
l_{31} & 0 & l_{33}%
\end{array}%
\right] ,\text{ }A_{6}^{T}=\left[ 
\begin{array}{ccc}
l_{11} & l_{21} & l_{31} \\ 
l_{12} & l_{22} & 0 \\ 
0 & 0 & l_{33}%
\end{array}%
\right] . 
\]
We consider the Lyapunov equation 
\begin{equation}
P_{5}A_{5}+A_{5}^{T}P_{5}=-I,\text{ }P_{5}=\left[ 
\begin{array}{ccc}
p_{11} & p_{12} & p_{13} \\ 
p_{21} & p_{22} & p_{23} \\ 
p_{31} & p_{32} & p_{33}%
\end{array}%
\right] ,\text{ }p_{ij}=p_{ji}.  \tag{6.22}
\end{equation}

\bigskip By solving $\left( 6.22\right) ,$ in a similar way as in the
Theorem 6.2 we obtain 
\[
p_{33}=-\frac{1}{2l_{33}},\text{ }p_{13}=\frac{d_{1}}{d},\text{ }p_{23}=%
\frac{d_{2}}{d},\text{ }p_{11}=\frac{D_{1}}{D},\text{ }p_{12}=\frac{D_{2}}{D}%
,\text{ }p_{22}=\frac{D_{3}}{D}, 
\]%
where 
\begin{equation}
\text{ }d_{1}=-\frac{l_{21}l_{31}}{2l_{33}},\text{ }d_{2}=\frac{l_{31}}{%
2l_{33}}\left( l_{11}+l_{33}\right) ,  \tag{6.23}
\end{equation}%
\[
\text{ }D_{1}=-\frac{1}{2}l_{21}^{2}+l_{22}\left( l_{11}+l_{22}\right)
\left( \frac{1}{2}+l_{31}p_{13}\right) + 
\]%
\[
\left( \frac{1}{2}+l_{31}p_{13}\right)
l_{12}l_{21}+l_{21}l_{22}l_{31}p_{23}, 
\]%
\[
\text{ }D_{2}=\frac{1}{2}l_{11}l_{21}+l_{12}l_{22}\left( \frac{1}{2}%
+l_{31}p_{13}\right) -l_{11}l_{22}l_{31}p_{23}, 
\]%
\[
\text{ }D_{3}=l_{11}l_{12}l_{31}p_{23}+\frac{1}{2}l_{12}l_{21}-\frac{1}{2}%
l_{11}\left( l_{11}+l_{22}\right) -k_{12}^{2}\left( \frac{1}{2}%
+l_{31}p_{13}\right) . 
\]

\bigskip Hence, the eigenvalues of $A_{6}$ are positive if the quadratic
function 
\[
V_{5}\left( x\right)
=X^{T}P_{2}X=p_{11}x_{1}^{2}+p_{22}x_{2}^{2}+p_{33}x_{3}^{2}+2p_{12}x_{1}x_{2}+ 
\]%
\[
2p_{13}x_{1}x_{3}+2p_{23}x_{2}x_{3} 
\]%
is positive defined. In a similar way we obtain that $V_{6}\left( x\right) $%
\ is positive defined, when $\frac{D_{k}}{D}>0,$ $k=1,3$ and the assumptions
of type $\left( 6.4\right) $\ are hold.

By reasoning as in the proof of Theorem 6.1 we obtain that the inequality 
\[
\dot{V}_{6}\left( x\right) =\dsum\limits_{k=1}^{3}\frac{\partial V_{6}}{%
\partial x_{k}}\frac{dx_{k}}{dt}\leq 0 
\]%
is valid if the following holds%
\[
p_{11}x_{1}+p_{12}x_{2}+p_{13}x_{3}\geq 0,\text{ }%
p_{12}x_{1}+p_{22}x_{2}+p_{23}x_{3}\geq 0,\text{ }%
p_{13}x_{1}+p_{23}x_{2}+p_{33}x_{3}\geq 0, 
\]%
\[
B_{1}\left( x_{1}\right) -D_{1}\left( x_{1},x_{2}\right) -h_{1}\left(
x_{1},x_{3}\right) \leq 0,\text{ }B_{2}\left( x_{2}\right) -D_{2}\left(
x_{2}\right) -h_{2}\left( x_{1},x_{2}\right) \leq 0, 
\]%
\begin{equation}
B_{3}\left( x_{1},x_{3}\right) -D_{3}\left( x_{3}\right) -h_{3}\left(
x_{1},x_{3}\right) \leq 0.  \tag{6.24}
\end{equation}

\bigskip \textbf{Remark 6.6. }By assumption $p_{33}>0$ and by $\left(
6.23\right) ,\ p_{kk}>0$ when $\frac{D_{k}}{D}>0,$ $k=1,3.$ Moreover, by
using $\left( 6.23\right) $ we can deduced the conditions on $l_{ij}$ that
the assumptions of type $\left( 6.4\right) $\ are hold.

\begin{center}
\textbf{7. Basins of multiphase attractions }
\end{center}

\bigskip In this section we will derived the domains of multipoint
attraction sets of the problem $\left( 1.3\right) -\left( 1.4\right) $ at
the the following attractor points $\left( 4.2\right) ,$ where $a_{\pm }$, $%
b_{\mp },$ $\bar{x}_{1},\bar{x},$ $x_{1i},x_{2j},x_{3ij}$ were defined by $%
\left( 4.16\right) $ and $\left( 4.24\right) .$

Lyapunov's method can be used to find the region of attraction or an
estimate of it. We show in this section the following results:

\textbf{Theorem 7.1. }Assume that the all conditions of Theorem 6.1 are
satisfied. Then the basin of multiphase attraction set of $\left( 1.3\right)
-\left( 1.4\right) $ at $\ \bar{x}=\left( 1,0,0\right) $ belongs to the set $%
\Omega _{C}\subset \Omega _{1},$ where $\Omega _{1}$ was defined by $\left(
4.8\right) $ and

\[
\Omega _{C}=\left\{ x\in R_{+}^{3}\text{: }V_{1}\left( x\right) \leq C\text{ 
}\right\} , 
\]%
here a positive constant $C$ is defined in bellow.

\textbf{Proof. }We are interested in the largest set $\Omega _{C}$ that we
can determine the largest value for the constant\textbf{\ }$C$\textbf{\ }%
such that $\Omega _{C}\subset D\left( V_{1}\right) ,$ where 
\[
D\left( V_{1}\right) =\left\{ x\in R^{3},\text{ }V_{1}\left( x\right) \geq 0,%
\text{ }\dot{V}_{1}\left( x\right) <0\right\} . 
\]%
\textbf{\ }Let us now, find the set $\Omega _{C}\subset B_{r}\left( \bar{x}%
\right) ,$ where 
\[
C<\min_{\left\vert x-\bar{x}\right\vert =r}V_{1}\left( x\right) =\lambda
_{\min }\left( P_{1}\right) r^{2}, 
\]%
here $P_{1}$ was defined by $\left( 4.1\right) $, $\lambda _{\min }\left(
P_{1}\right) $ denotes a minimum eigenvalue of the corresponding matrix $%
A_{1}$.

Moreover, for some $C>0$ the inclusion$\ \Omega _{C}\subset \Omega _{1}$
means the existence of $C>0$ such that $x\in \Omega _{C}$ implies $x\in
\Omega _{1}$, where 
\[
\text{ }\Omega _{1}=\left\{ x\in \mathbb{R}_{+}^{3}\text{, }%
x_{j}=x_{j_{0}}+\dsum\limits_{k=1}^{m}\alpha _{jk}x_{j}\left( t_{k}\right)
\geq 0\text{, }j=1,2,3,\text{ }x_{2}\geq \eta _{2},\right. 
\]%
\[
p_{11}x_{1}+p_{12}x_{2}+p_{13}x_{3}\geq 0,\text{ }%
p_{12}x_{1}+p_{22}x_{2}+p_{23}x_{3}\geq 0,\text{ }%
p_{13}x_{1}+p_{23}x_{2}+p_{33}x_{3}\geq 0, 
\]%
\[
B_{1}\left( x_{1}\right) -D_{1}\left( x_{1},x_{2}\right) -h_{1}\left(
x_{1},x_{3}\right) \leq 0,\text{ }B_{2}\left( x_{2}\right) -D_{2}\left(
x_{2}\right) -h_{2}\left( x_{1},x_{2}\right) \leq 0, 
\]%
\begin{equation}
\left. B_{3}\left( x_{1},x_{3}\right) -D_{3}\left( x_{3}\right) -h_{3}\left(
x_{1},x_{3}\right) \leq 0\right\} .  \tag{7.1}
\end{equation}%
here $O_{\delta }\left( t_{0}\right) $ was defined by $\left( 1.3\right) $, $%
p_{ij}$, $a_{ij}$ were defined by $\left( 6.2\right) $ and $\left(
4.8\right) ,$ respectively,\ i.e. 
\[
\text{ }p_{33}=-\frac{1}{2a_{33}}\text{, }p_{13}=\frac{a_{31}}{2\left(
a_{11}+a_{33}\right) a_{33}},\text{ }p_{11}=-\frac{1}{a_{11}}\left( \frac{1}{%
2}+a_{31}p_{13}\right) ,\text{ } 
\]

\[
p_{23}=-\frac{a_{12}p_{13}}{a_{22}+a_{33}},\text{ }p_{12}=-\frac{\left(
a_{12}p_{11}+a_{31}p_{23}\right) }{\left( a_{11}+a_{22}\right) }\text{, }%
p_{22}=-\frac{-\left( \frac{1}{2}+a_{12}p_{12}\right) }{a_{22}}. 
\]
\[
a_{11}=\frac{\partial }{\partial x_{1}}\left[ B_{1}-D_{1}\right] \left(
0\right) -\frac{\partial h_{1}}{\partial x_{1}}\left( 0\right) \text{, }%
a_{12}=-\frac{\partial D_{1}}{\partial x_{1}}\left( 0\right) , 
\]%
\[
\text{ }a_{22}=\frac{d}{dx_{2}}\left[ B_{2}-D_{2}\right] \left( 0\right) ,%
\text{ }a_{31}=\frac{\partial B_{3}}{\partial x_{1}}\left( 0\right) -\frac{%
\partial h_{1}}{\partial x_{1}}\left( 0\right) \text{,} 
\]
\[
a_{33}=\frac{d}{dx_{3}}\left[ B_{3}-D_{3}\right] \left( 0\right) . 
\]

\textbf{Remark 7.1.} By assumptions of theorem \ $p_{ii}>0.$ By Remark 6.1
if $a_{31}>0,$ then $p_{13}>0;$ moreover, $p_{23}>0,$ $p_{12}>0$ when $%
a_{31}>0$ and $a_{12}>0.$ Then $\left( 7.1\right) $ holds if 
\[
B_{1}\left( x_{1}\right) -D_{1}\left( x_{1},x_{2}\right) -h_{1}\left(
x_{1},x_{3}\right) \leq 0,\text{ }B_{2}\left( x_{2}\right) -D_{2}\left(
x_{2}\right) -h_{2}\left( x_{1},x_{2}\right) \leq 0, 
\]%
\begin{equation}
\left. B_{3}\left( x_{1},x_{3}\right) -D_{3}\left( x_{3}\right) -h_{3}\left(
x_{1},x_{3}\right) \leq 0\right\} .  \tag{7.2}
\end{equation}

\bigskip In view of $\left( 4.8\right) ,$ $a_{31}>0$ , $a_{12}>0,$ when $%
\frac{\partial B_{3}}{\partial x_{1}}\left( 0\right) >\frac{\partial h_{1}}{%
\partial x_{1}}\left( 0\right) $ and $\frac{\partial D_{1}}{\partial x_{1}}%
\left( 0\right) <0.$

\bigskip

\ Hence, \ 
\[
\text{ }\Omega _{10}=\left\{ x\in \mathbb{R}_{+}^{3},\right. \left.
b_{11}\left( x_{1}-1\right) ^{2}+\left( b_{22}+b_{12}\right)
x_{2}^{2}+x_{3}^{2}\leq \right\} \text{ } 
\]%
\[
b_{11}+\left( \beta _{1}+\beta _{2}\eta _{2}\right) ^{2}\text{, }x_{1}\geq
1\left. {}\right\} \subset \Omega _{1}. 
\]%
So, it is not hard to see that 
\[
B_{\bar{r}}\left( \bar{x}\right) =\left\{ x\in R^{3}\text{, }\left\vert x-%
\bar{x}\right\vert <\bar{r}\right\} \subset \Omega _{1}, 
\]%
where 
\[
\text{ }\tilde{r}=\eta _{0}^{\frac{1}{2}}\left[ b_{11}+\left( \beta
_{1}+\beta _{2}\eta _{2}\right) ^{2}\right] ^{\frac{1}{2}},\text{ }\eta
_{0}=\max \left\{ b_{11},\text{ }b_{22}+b_{12},1\right\} . 
\]%
Then we obtain 
\[
C<\min_{\left\vert x\right\vert =r_{1}}V_{1}\left( x\right) =\lambda _{\min
}\left( P_{1}\right) \tilde{r}^{2}, 
\]%
i.e. 
\[
C<\lambda _{\min }\left( P_{1}\right) r_{0}^{2},\text{ }r_{0}=\min \left\{ r,%
\text{ }\tilde{r}\right\} . 
\]

Now, we consider the equilibria point $E_{2}\left( 0,1,0\right) $ and prove
the following result

\textbf{Theorem 5.2. }Assume that the all conditions of Theorem 4.2 and $%
\left( 4.15\right) $\ are satisfied.Then the basin of multiphase attraction
set of $\left( 1.3\right) -\left( 1.4\right) $ at $E_{2}\left( 0,1,0\right) $
is whole $\mathbb{R}_{+}^{3}.$

\textbf{Proof. }Indeed, by Theorem 4.2 the system $\left( 1.3\right) $ is
global stabile at $E_{2}\left( 0,1,0\right) .$ Thus, the basin of multiphase
attraction set coincides with $\mathbb{R}_{+}^{3}.$

\textbf{Theorem 5.3. }Assume that the all conditions of Theorem 4.3 are
satisfied. Then the basin of multiphase attraction set of $\left( 1.3\right)
-\left( 1.4\right) $ at $E_{3}\left( a_{\pm },0,b_{\mp }\right) $ belongs to
the set $\Omega _{C}\subset \Omega _{3},$ where $\Omega _{3}$ was defined by 
$\left( 4.23\right) ,$ here $V_{3}\left( x\right) $ was defined by $\left(
4.15\right) .$

\textbf{Proof. }We will find $C>0$ such that $\Omega _{C}\subset B_{r}\left(
E_{3}\right) \cap \Omega _{3}$. It is clear to see that $\Omega _{C}\subset
B_{r}\left( E_{3}\right) $ for 
\[
C<\min_{\left\vert x-\bar{x}\right\vert =r}V_{3}\left( x\right) =\lambda
_{\min }\left( P_{3}\right) r^{2},\text{ }\bar{x}=\left( a_{\pm },0,b_{\mp
}\right) , 
\]%
here $\lambda _{\min }\left( P_{3}\right) $ denotes a minimum eignevalue of $%
A_{3}$. Let $\Omega _{3}$ is a domain defined by $\left( 4.23\right) $, i.e. 
\[
\Omega _{3}=\left\{ x\in \mathbb{R}_{+}^{3}\text{: }x_{j}=x_{j_{0}}+\dsum%
\limits_{k=1}^{m}\alpha _{jk}x_{j}\left( t_{k}\right) \geq 0\text{, }%
j=1,2,3,\right. 
\]%
\[
\alpha _{1}x_{1}+\alpha _{2}x_{2}+\alpha _{3}x_{3}\geq \gamma _{0}\text{, }%
x_{1}\geq \gamma _{1},\text{ }x_{2}\leq \gamma _{2}x_{3},\text{ }x_{3}\leq
\gamma _{3}x_{1}, 
\]%
\[
\left( b_{11}+b_{11}a_{\pm }+b_{13}b_{\mp }\right) \left( x_{1}-a_{\pm
}\right) ^{2}+r_{2}\left( b_{12}a_{\pm }+b_{23}b_{\mp }+b_{22}\right)
x_{2}^{2}\leq 
\]%
\[
\left. r_{2}\left( b_{12}a_{\pm }+b_{23}b_{\mp }+b_{22}\right)
+b_{11}x_{1}^{3}\text{,}\right\} , 
\]%
where%
\[
\alpha _{1}=\min \left\{ {}\right. \left[ b_{11}a_{\pm }+b_{13}b_{\mp
}-2a_{\pm }\left( b_{11}+b_{11}a_{\pm }+b_{13}b_{\mp }\right) \right] \text{%
, } 
\]%
\[
\text{ }\left. b_{11}a_{12}+b_{12}a_{21},\text{ }b_{12},\text{ }%
b_{13}\right\} , 
\]%
\[
\alpha _{2}=\text{ }\min \left\{ {}\right. r_{2}\left( b_{12}a_{\pm
}+b_{23}b_{\mp }\right) -2a_{\pm }r_{2}\left( b_{12}a_{\pm }+b_{23}b_{\mp
}+b_{22}\right) \text{,} 
\]%
\[
\left. b_{12}a_{12}+b_{22}a_{21},\text{\ }b_{22},\text{ }b_{23}\right\} ,%
\text{ }\alpha _{3}=\min \left\{ b_{13}a_{12},\text{ }b_{23},\text{ }%
b_{33}\right\} =b_{23}, 
\]%
\[
\gamma _{0}=\left( b_{11}a_{12}a_{\pm }+b_{12}+b_{13}a_{12}b_{\mp
}+a_{21}b_{12}a_{\pm }+a_{21}b_{23}b_{\mp }\right) , 
\]%
\[
\gamma _{1}=\frac{\left( b_{11}a_{\pm }+b_{13}b_{\mp }\right) a_{13}+b_{13}}{%
\left( b_{13}+a_{13}b_{11}\right) },\text{ }\gamma _{2}=\frac{%
a_{13}b_{13}x_{3}}{-a_{21}b_{23}}\text{, }\gamma _{3}=\frac{b_{11}a_{13}}{%
-b_{23}a_{21}}. 
\]%
\ It is clear that $\alpha _{2}$, $\alpha _{3}\leq 0$ and $\alpha _{1}>0.$
Hence, $\alpha _{1}x_{1}-\gamma _{0}>0.$\ Moreover, since 
\[
\alpha _{1}x_{1}+\alpha _{2}x_{2}+\alpha _{3}x_{3}\geq \gamma _{0}\text{, }%
x_{1}\geq \gamma _{1},\text{ }x_{2}\leq \gamma _{2}x_{3},\text{ }x_{3}\leq
\gamma _{3}x_{1} 
\]%
we get \ 
\[
0\leq x_{3}\leq \beta _{1}\gamma _{1}-\beta _{2}, 
\]%
where 
\[
\beta _{1}=\frac{\alpha _{1}}{-\left( \alpha _{2}\gamma _{2}+\alpha
_{3}\right) }\text{, }\beta _{2}=\frac{\gamma _{0}}{-\left( \alpha
_{2}\gamma _{2}+\alpha _{3}\right) }. 
\]%
Thus,%
\begin{equation}
\Omega _{30}=\left\{ x\in \mathbb{R}_{+}^{3}\text{:}\right. \text{ }%
x_{j}=x_{j_{0}}+\dsum\limits_{k=1}^{m}\alpha _{jk}x_{j}\left( t_{k}\right)
\geq 0\text{, }j=1,2,3,  \tag{5.3}
\end{equation}%
\[
\left( b_{11}+b_{11}a_{\pm }+b_{13}b_{\mp }\right) \left( x_{1}-a_{\pm
}\right) ^{2}+r_{2}\left( b_{12}a_{\pm }+b_{23}b_{\mp }+b_{22}\right)
x_{2}^{2}+x_{3}^{2}\leq 
\]%
\[
\left. r_{2}\left( b_{12}a_{\pm }+b_{23}b_{\mp }+b_{22}\right) +b_{11}\gamma
_{1}^{3}+\left( \beta _{1}\gamma _{1}-\beta _{2}\right) ^{2}\right\} . 
\]%
From $\left( 4.23\right) $ it is not hard to see that%
\[
B_{\bar{r}}\left( \bar{x}\right) =\left\{ x\in R_{+}^{3}\text{, }\left\vert
x-\bar{x}\right\vert <\bar{r}\right\} \subset \Omega _{3}\text{ for }\bar{x}%
=\left( 0,a_{\pm },b_{\mp }\right) , 
\]%
where 
\[
\left( \bar{r}\right) ^{2}=\frac{1}{\eta }\left[ r_{2}\left( b_{12}a_{\pm
}+b_{23}b_{\mp }+b_{22}\right) +b_{11}\gamma _{1}^{3}+\left( \beta
_{1}\gamma _{1}-\beta _{2}\right) ^{2}\right] , 
\]%
\[
\eta =\max \left\{ \left( b_{11}+b_{11}a_{\pm }+b_{13}b_{\mp }\right) ,\text{
}r_{2}\left( b_{12}a_{\pm }+b_{23}b_{\mp }+b_{22}\right) \text{, }1\right\}
. 
\]%
Then we obtain that 
\[
C<\min_{\left\vert x-\bar{x}\right\vert =\bar{r}}V_{3}\left( x\right)
=\lambda _{\min }\left( P_{3}\right) \bar{r}^{2}, 
\]%
i.e. 
\[
C<\lambda _{\min }\left( P_{3}\right) \bar{r}^{2}\text{ for }r_{0}=\min
\left\{ r,\text{ }\bar{r}\right\} . 
\]%
Consider the point $E_{4}\left( \bar{x}_{1},\bar{x}_{2},0\right) .$ By
reasoning as the above we prove the following result:

\textbf{Theorem 5.4. }Assume that the all conditions of Theorem 4.4 are
satisfied. Then the basin of multiphase attraction sets of $\left(
1.3\right) -\left( 1.4\right) \ $at\ $E_{4}\left( \bar{x}_{1},\bar{x}%
_{2},0\right) $ belongs to the set $\Omega _{4},$ where $\Omega _{4}$ was
defined by $\left( 4.31\right) .$

\textbf{Proof. }We will find $C>0$ such that $\Omega _{C}\subset B_{r}\left(
E_{4}\right) \subset \Omega _{4}$. It is clear to see that $\Omega
_{C}\subset B_{r}\left( \text{ }\bar{x}\right) $ for 
\[
C<\min_{\left\vert x-\bar{x}\right\vert =r}V_{4}\left( x\right) =\lambda
_{\min }\left( P_{4}\right) r^{2},\text{ }\bar{x}=\left( \bar{x}_{1},\bar{x}%
_{2},0\right) , 
\]%
here $\lambda _{\min }\left( P_{4}\right) $ denotes a minimum eigenvalue of $%
A_{4}.$ From $\left( 4.31\right) $ we get 
\begin{equation}
\Omega _{40}=\left\{ {}\right. x\in \mathbb{R}_{+}^{3}\text{: }%
x_{j}=x_{j_{0}}+\dsum\limits_{k=1}^{m}\alpha _{jk}x_{j}\left( t_{k}\right)
\geq 0\text{, }j=1,2,3,\text{ }  \tag{5.4}
\end{equation}%
\[
x_{1}\leq \gamma _{1},\text{ }x_{2}\geq \gamma _{2},\text{ }x_{3}\leq \gamma
_{3}, 
\]%
\[
\left( b_{11}\bar{x}_{1}+b_{12}\bar{x}_{2}\right) \left( x_{1}-\bar{x}%
_{1}\right) ^{2}+r_{2}\left( b_{12}\bar{x}_{1}+b_{22}\bar{x}_{2}\right)
\left( x_{2}-\bar{x}_{2}\right) ^{2}\leq 
\]%
\[
\left( b_{11}\bar{x}_{1}+b_{12}\bar{x}_{2}\right) \bar{x}_{1}^{2}+r_{2}%
\left( b_{12}\bar{x}_{1}+b_{22}\bar{x}_{2}\right) \bar{x}%
_{2}^{2}+b_{22}r_{2}x_{2}^{3}\text{, }x_{3}\leq \frac{a_{21}b_{23}}{%
-b_{13}a_{13}}x_{2}, 
\]%
\[
\left. \alpha _{1}x_{1}+\alpha _{2}x_{2}+\alpha _{3}x_{3}\geq b_{13}\right\}
\subset \Omega _{4}, 
\]%
where 
\[
\gamma _{1}=\text{\ }\frac{b_{12}-r_{2}\left( b_{12}\bar{x}_{1}+b_{22}\bar{x}%
_{2}\right) }{b_{12}},\text{ }\gamma _{3}=\frac{\left( b_{11}\bar{x}%
_{1}+b_{12}\bar{x}_{2}\right) }{a_{13}\left( b_{11}\bar{x}_{1}+b_{12}\bar{x}%
_{2}\right) }, 
\]%
\[
\text{ }\gamma _{2}=\max \left\{ \frac{a_{21}\left( b_{12}\bar{x}_{1}+b_{22}%
\bar{x}_{2}+b_{12}r_{2}\right) }{\left( a_{12}b_{12}+a_{21}b_{22}\right) }%
\text{, }1\text{, }\frac{a_{12}\left( b_{11}\bar{x}_{1}+b_{12}\bar{x}%
_{2}\right) }{b_{12}}\right\} , 
\]%
\[
\alpha _{1}=\text{ }\min \left\{ b_{11},\text{ }b_{13}\text{ }\right\} ,%
\text{ }\alpha _{3}=\text{ }\min \left\{ b_{13},\text{ }b_{23}\right\} , 
\]%
\[
\alpha _{2}=\min \left\{ \left( b_{12}+a_{12}b_{11}+a_{21}b_{22}\right) 
\text{, }a_{12}b_{12}+a_{21}b_{22}\text{, }b_{23}\right\} . 
\]

From $\left( 5.4\right) $ It is not hard to see that $\gamma _{1}\leq \frac{%
\alpha _{2}\gamma _{2}}{-b_{13}}$ and

\[
B_{r}\left( \bar{x}\right) =\left\{ x\in R_{+}^{3}\text{, }\left\vert x-\bar{%
x}\right\vert <\bar{r}\right\} \subset \Omega _{40}\text{ for }\bar{x}%
=\left( \bar{x}_{1},\bar{x}_{2},0\right) , 
\]%
where 
\[
\left( \bar{r}\right) ^{2}=\frac{1}{\eta }\left[ \left( b_{11}\bar{x}%
_{1}+b_{12}\bar{x}_{2}\right) \bar{x}_{1}^{2}+r_{2}\left( b_{12}\bar{x}%
_{1}+b_{22}\bar{x}_{2}\right) \bar{x}_{2}^{2}+b_{22}r_{2}\gamma
_{2}^{3}+d^{2}\right] , 
\]%
\[
\eta =\max \left\{ b_{11}\bar{x}_{1}+b_{12}\bar{x}_{2},\text{ }r_{2}\left(
b_{12}\bar{x}_{1}+b_{22}\bar{x}_{2}\right) \text{, }1\right\} \text{, }%
d=\min \left\{ \frac{\alpha _{2}\gamma _{2}}{-b_{13}}-\gamma _{1},\text{ }%
\gamma _{3}\right\} . 
\]%
Then we obtain that 
\[
C<\min_{\left\vert x-\bar{x}\right\vert =\bar{r}}V_{4}\left( x\right)
=\lambda _{\min }\left( P_{4}\right) \bar{r}^{2}, 
\]%
i.e. 
\[
C<\lambda _{\min }\left( P_{4}\right) \bar{r}^{2}\text{ for }r_{0}=\min
\left\{ r,\text{ }\bar{r}\right\} . 
\]

Consider the points $E_{ij}.$

\textbf{Theorem 5.5. }Assume that the all conditions of Theorem 4.5 are
satisfied. Then the basin of multiphase attraction sets of $\left(
1.3\right) -\left( 1.4\right) $ at points $E_{ij}$ belong to the $\Omega
_{ij},$ where $\Omega _{ij}$ was defined by $\left( 4.38\right) .$

\textbf{Proof. }We will find $C>0$ such that $\Omega _{C}\subset B_{r}\left(
E_{ij}\right) \subset \Omega _{ij}$. It is clear to see that $\Omega
_{C}\subset B_{r}\left( \text{ }\bar{x}\right) $ for 
\[
C<\min_{\left\vert x-\bar{x}\right\vert =r}V_{5}\left( x\right) =\lambda
_{\min }\left( P_{5}\right) r^{2}, 
\]%
here $\lambda _{\min }\left( P_{5}\right) $ denotes a minimum eignevalue of $%
A_{5}.$ Assume $a_{13}>1.$Then from $\left( 4.38\right) $ it is not hard to
see that

\begin{equation}
B_{r}\left( E_{ij}\right) \subset \Omega _{ij0}=\left\{ x\in \mathbb{R}%
_{+}^{3}\text{: }x_{j}=x_{j_{0}}+\dsum\limits_{k=1}^{m}\alpha
_{jk}x_{j}\left( t_{k}\right) \geq 0\text{, }j=1,2,3,\right.  \tag{5.5}
\end{equation}%
\[
x_{1}\leq \gamma _{1},\text{ }x_{2}\geq 1,\text{ }x_{3}\leq \frac{1}{a_{13}},%
\text{ } 
\]%
\[
\text{\ }Q_{1}\left( x_{1}-x_{1i}\right) ^{2}+Q_{2}r_{2}\left(
x_{2}-x_{2j}\right) ^{2}+\left( x_{3}-x_{3ij}\right) ^{2}\leq
Q_{1}x_{1i}^{2}+Q_{1}x_{2j}^{2} 
\]%
\[
\text{\ }\left. \text{\ }+\left( \frac{1}{a_{13}}-x_{3ij}\right)
^{2}+p_{22}r_{2}+d^{2},\text{ }-\text{ }\left[ \alpha _{1}x_{1}+\alpha
_{2}x_{2}\right] \leq \alpha _{3}x_{3}\right\} ,\text{ } 
\]%
where%
\[
\alpha _{1}=\text{ }\min \left\{ p_{11},\text{ }p_{23}a_{21}+p_{13}a_{13},%
\text{ }p_{12}a_{21}\text{, }p_{13}\text{ }\right\} , 
\]%
\[
\text{ }\alpha _{2}=\text{ }\min \left\{ p_{11}a_{12}+p_{12}\text{, }%
p_{12}a_{13}\text{, }p_{12}\left( a_{12}+r_{2}\right) +p_{22}a_{21}\text{, }%
p_{23}\right\} , 
\]%
\[
\alpha _{3}=\text{ }\min \left\{ p_{11}a_{13},\text{ }p_{13}a_{13},\text{ }%
p_{13}a_{12}\text{, }p_{33}\right\} ,\text{ }d=\frac{-p_{12}}{\alpha _{3}}%
\left( 1+\gamma _{1}\right) \text{, } 
\]%
\[
a=\max \left\{ a_{21},\text{ }a_{12}r_{2}\right\} 
\]%
\[
\gamma _{1}=\frac{r_{2}}{\left( a_{13}+2x_{1i}\right) Q_{1}+\left(
a_{21}+r_{2}+2x_{2j}\right) Q_{2}}, 
\]

\[
\left( \bar{r}\right) ^{2}=\frac{1}{\eta }\left[
Q_{1}x_{1i}^{2}+Q_{1}x_{2j}^{2}+\left( \frac{1}{a_{13}}-x_{3ij}\right)
^{2}+p_{22}r_{2}+d^{2}\right] \text{, } 
\]%
\[
\eta =\max \left\{ Q_{1},\text{ }Q_{2},\text{ }1\right\} . 
\]%
Then we obtain that 
\[
C<\min_{\left\vert x-\bar{x}\right\vert =\bar{r}}V_{5}\left( x\right)
=\lambda _{\min }\left( P_{5}\right) \bar{r}^{2}, 
\]%
i.e. 
\[
C<\lambda _{\min }\left( P_{5}\right) \bar{r}^{2}\text{ for }r_{0}=\min
\left\{ r,\text{ }\bar{r}\right\} . 
\]

\bigskip \textbf{Conclusion.} Taking into account different and effective
features of mathematical modelling and its possibilities to figure out a
problem in dynamics on the basis of its logic properties, it was surely
pointed out the characteristics of a mathematical model to use in
description of needed processes of a given dynamic system with identified
problems. In this paper, a three dimensional model was devoted to
mathematical description and regulation possibilities of uncontrolled tumor
processes by organism as a complex system. The dynamics of interactions of
the dimensions corresponded to tumor cells, immune cells and healthy --
\textquotedblleft host\textquotedblright\ -- cells were given as forces of
vectors, negatively or positively converging to basins of attractions,
depending on their importance for the complex system. In order to make the
model subjected to control, there was included multiphase IVP, describing
the system's important parameters to operate with it in the farther
processes of stages of development. The model was undergone different
changes to determine its limits of survival: it was determined the
conditions of boundedness the system can be restricted, invariance in non-
negativity, which means the model keeps its properties of reactions to
changing in proper way, being subjected to different analysis, and the
circumstances the system can be forced to be dissipated in. The system was
exposed to changing pressures to estimate its convenience to biologically
important properties as points of equilibria and Lyapunov stability
conditions. The next step in exploring of the model were very complex and
logistic approaches to its properties for verification of the conditions,
providing the global equilibria points and multimodal attraction sets,
having biologically strong value in regulation of the processes towards the
positive effects of feasible medical external implementation at the
convenient stages, determined by multimodal attraction basins.

\textbf{Biological implications. }Here we study a multiphase host-tumor
model that enhances the type of effector immune cells that can fight a
tumor, and stimulates effector immune cells to proliferate.\textbf{\ }%
Interactions between cancer tumor cells, healthy host cells and the effector
immune cells can explain long-term tumor relapse. Here, the sufficient
conditions is derived that under which the possible biologically feasible
dynamics is stable in the Lyapunov sense, and a converges to one of
equilibrium points. Since these equilibrium points have a biological sense,
we notice that understanding limit properties of dynamics of cells
populations based on solving the problem $(1.3)-\left( 1.4\right) $ may be
of an essential interest for the prediction of health conditions of a
patient without a treatment, when the data (e.g. the status of blood cells
shown above) that determines the condition of the patient are compared at
various times $t_{0},t_{1},...,t_{m}$ and correlated. In the section 3, we
find the positively invariant domain $B_{\alpha ,m}$ that depend on
multipoint IVP condition parameters $\alpha _{k}$, $t_{k}$ and $m.$
Moreover, the boundedness of orbits of the system $\left( 1.3\right) -\left(
1.4\right) $ is derived. As a result, the future evolution of cells
populations involved in this model is completely predictable in the
following sense: by knowing the specific linear connection between the
tumor, guest and immune cells at the $t_{0},$ $t_{1},$...$t_{m}$ time phase
densities, populations has an accurate and predictable estimate of its
change. In the section 4, lyapunov stability of the system $\left(
1.3\right) $ at the corresponding equilibria points are studied. We show
that the system $\left( 1.3\right) $ is global stable at the "free tumor "
equilibria\ point $E_{2}\left( 0,1,0\right) .$ In the section 5, the basins
of multiphase attractors of the system $\left( 1.3\right) -\left( 1.4\right) 
$ (dependent on multipoint parameters of IVP) are constructed.

\textbf{Acknowledgements}

The author is thanking to Assist. Prof. Department of Biophysics Yeditepe
University A. Maharramov, Assoc. Prof. Department of Immunology Yeditepe
University G. Yanikkaya Demirel, Assoc. Prof. Department of Medical
Microbiology Yeditepe University \.{I}brahim Ch. Acuner and Prof. Dr.
Faculty of Health Sciences Okan University Aida Sahmurova according to their
valuable suggestions in the field of medicine and biology.

\bigskip \textbf{References}

\begin{enumerate}
\item Kuznetsov V. A., Makalkin I. A., Taylor M. A., Perelson A. S.,
Nonlinear dynamics of immunogenic tumors: parameter estimation and global
bifurcation analysis, Bull. Math. Biol. 1994(56), 295--321.

\item Adam J. A, Bellomo C., A survey of models for tumor-immune system
dynamics, Boston, MA: Birkhauser, 1996.

\item Eftimie R, Bramson J. L., Earn D. J. D., Interactions between the
immune system and cancer: a brief review of non-spatial mathematical models,
Bull. Math. Biol. 2011(73), 2--32.

\item Kirschner D, Panetta J., Modelling immunotherapy of the tumor--immune
interaction, J. Math. Biol. 1998(37), 235--52.

\item de Pillis L. G., Radunskaya A., The dynamics of an optimally
controlled tumor model: a case study, Math. Comput. Modell 2003(37),1221--44.

\item Nani F., Freedman H. A., Mathematical model of cancer treatment by
immunotherapy, Math Biosci. 2000(163), 159--99.

\item Owen M. and Sherratt J., Mathematical modelling macrophage dynamics in
tumors, Mathematical Models and Methods in Applied Sciences 9 (4)(1999),
513-539.

\item Chaplain M.A.J., Special issue on mathematical models for the growth,
development and treatment of tumours, Math. Models Meth. Appt. Sci. 9 (1999).

\item Arciero, J., Jackson, T., \& Kirschner, D., A mathematical model of
tumor-immune evasion and siRNA treatment, Discrete Contin. Dyn. Syst., Ser.
B, 4(1) (2004), 39--58.

\item Kirschner D, Tsygvintsev A., On the global dynamics of a model for
tumor immunotherapy, Mathematical Biosciences and engineering, (6)3 (2009),
573-583.

\item Starkov, K. E., Krishchenko, A. P., On the global dynamics of one
cancer tumour growth model, Commun. Nonlinear Sci. Numer. Simul. 19 (2014),
1486--1495.

\item Itik I. M., Banks S. P., Chaos in a three-dimensional cancer model,
Int J. Bifurcation Chaos 2010, 2010(20), 71--79.

\item Levine, H., A., Pamuk S., Sleeman B. D., Mathematical modeling of
capillary formation and development in tumor angiogenesis, Penetration into
the Stroma bulletin of Mathematical Biology (2001) 63, 801--863.

\item Jackson T., Komarova N. and Swanson K., Mathematical oncology: Using
mathematics to enable cancer discoveries, American Mathematical Monthly,
121(9), (2014), 840-856.

\item Firmani B., Guerri L. and Preziosi L., Tumor/immune system competition
with medically induced activation/deactivation, Mathematical Models and
Methods in Applied Sciences 4 (9)(1999), 491-512 .

\item Gallas M. R., Gallas Marcia R. and Gallas J. A.C., Distribution of
chaos and periodic spikes in a three-cell population model of cancer, Eur.
Phys. J. Special Topics, 223 (2014), 2131--2144.

\item Itik I. M., Salamci M.U., Banks S. P., Optimal control of drug therapy
in cancer treatment, Nonlinear Analysis, 71 (2009), 1473--86.

\item El-Gohary A., Chaos and optimal control of equilibrium states of tumor
system with drug, Chaos, Solitons and Fractals, 41(2019), 425--435.

\item Iarosz, K. C , Borges, F. S., Batista, A. M., Baptista, M. S ,
Siqueira, R. A. N , Viana, R. L., Lopes, S. R \& Baptista M. D. S,
Mathematical model of brain tumour with glia-neuron interactions and
chemotherapy treatment, Journal of Theoretical Biology, 368(2015), 113-121.

\item Kuang Y., Nagy J. D. and Eikenberry S. E., Introduction to
mathematical oncology, Chapman \& Hall/CRC Mathematical and computational
biology series, CRC Press, Boca Raton (2016).

\item Bellomo N. and Preziosi L., Modelling and mathematical problems
related to tumor evolution and Its Interaction with the Immune system,
Mathematical and Computer Modelling 32 (2000) 413452.

\item Khalil H. Nonlinear systems, NJ: Prentice Hall, 2002.

\item Verhulst, F., Nonlinear differential equations and dynamical systems,
Springer-Verlag Berlin Heidelberg 1996.

\item Carl S., Heikkila S., Fixed point theory in ordered sets and
applications, Springer, 2010.
\end{enumerate}

\bigskip

\end{document}